\documentclass[12pt]{article} 

\usepackage[francais,english]{babel} 

\catcode`\ =\active \def {\`e} \catcode`\ˆ =\active \def ˆ{\`a} \catcode`\ =\active \def {\`u}
\catcode`\Ž =\active \def Ž{\'e} \catcode`\ =\active \def {\c c} \def \parn{\par\noindent} 
\catcode`\ =\active \def {\^e} \catcode`\™ =\active \def ™{\^o} 
\catcode`\" =\active \def "{\^\i } 

\catcode`\: =\active \def :{\ifmmode\string:\else\unskip\kern 2pt\string: \ignorespaces\fi } 
\catcode`\; =\active \def ;{\ifmmode\string;\else\unskip\kern 2pt\string; \ignorespaces\fi } 

\font \db = msbm10 at 12 pt

            \def \G{\Gamma}
   \def \s{\sigma}     
\def \e{\varepsilon}  \def \f{\varphi}    \def \rr{\varrho}  

\def \R{\mbox{\db R}}    \def \P{\mbox{\db P}} 
\def \N{\mbox{\db N}}

      \def \Z{\mbox{\db Z}}

\def \AA{{\cal A}}  \def \BB{{\cal B}}  \def \CC{{\cal C}}       \def \EE{{\cal E}}
           
\def \LL{{\cal L}}     \def \NN{{\cal N}}    \def \O{{\cal O}}

       \def \ol{\overline} 

\def \p{\partial}        \def \ii{\infty}  \def \sm{\setminus}  
  \def \lra{\longrightarrow}   \def \ub{\underbar}

   \def \ds{\displaystyle}   \def \ts{\textstyle}       \def \ss{\scriptstyle}
\def \rt1{\sqrt{-1}\,\,}  \def \1{^{-1}}            \def \2{^{-2}}             \def \5{{\ts {1\over 2}}}
\def \moins{{\ss \sm}}    \def\indf{\leavevmode\indent }

\def\cotg{{\rm cotg}\,}   \def\tg{{\rm tg}\,}   \def\coth{{\rm coth}\,}        
\def\ch{{\rm ch}\,}       \def\sh{{\rm sh}\,}   \def\arctg{{\rm Arctg}\,}   \def\argsh{{\rm argsh}\,}

\begin{document}

\newtheorem{defi}{Definition}
\newtheorem{theo}{Theorem}
\newtheorem{prop}{Proposition}  \newtheorem{propr}{Property}
\newtheorem{cor}{Corollary}
\newtheorem{lem}{Lemma}
\newtheorem{rem}{Remark} 

\newtheorem{formu}{} 
\newcounter{form}{} 

\newcommand \beq{\begin{equation}} \newcommand \eeq{\end{equation}} 

\newcommand \bthe{\begin{theo}}  \newcommand \ethe{\end{theo}}   
\newcommand \bpro{\begin{prop}}  \newcommand \epro{\end{prop}}   
\newcommand \bcor{\begin{cor}}    \newcommand \ecor{\end{cor}}     
\newcommand \blem{\begin{lem}}   \newcommand \elem{\end{lem}}   
\newcommand \brem{\begin{rem}}  \newcommand \erem{\end{rem}} 
\newcommand \bdefi{\begin{defi}}   \newcommand \edefi{\end{defi}}

\title{\bf  Relativistic Diffusion in G\"odel's Universe}
\author{Jacques FRANCHI}
\date { February 2007 }
\maketitle 

\vspace{-9mm} 

\begin{abstract}  K. G\"odel [G1] discovered his celebrated solution to Einstein equations in 1949. Additional contributions were made by Kundt [K] and Hawking-Ellis ([H-E], 5.7).  On the other hand, a  general Lorentz invariant operator, associated to the so-called ``relativistic diffusion'', and making sense in any Lorentz manifold, was introduced by Franchi-Le Jan in [F-LJ]. \   Here is proposed a first study of the relativistic diffusion in the framework of G\" odel's universe, which contains matter. 
\end{abstract}
 
\vspace{-5mm} 
   
\section{Introduction} \label{sec.I} \indf 
   K. G\"odel [G1] published his celebrated exact solution to Einstein equations in 1949. The most striking feature of this cosmological model was to be non-causal (though locally and geodesically causal), containing closed timelike curves. For this reason, it is generally considered as rather unphysical. Possessing a series of interesting properties, it aroused however a great interest among physicists. For example, it contains rotating matter, but no singularity. Moreover, the explicit exact solutions to Einstein equations are not so many. \par 
   W. Kundt [K] studied its geodesics, and S. Hawking and G. Ellis ([H-E], section 5.7) stress on coordinates (defined by G\"odel himself) showing up its rotational symmetry (about any point), to draw a nice picture of its dynamics. D. Malament ([M1],[M2]) calculated the minimal energy of a closed timelike curve. K. G\"odel [G2] discussed other rotating universes, which are spatially homogeneous, finite, and expanding, and he showed in particular that there exist a lot of strongly causal such cosmological models. 
   
    \par   \smallskip 
    The purpose of the present work is to study, in the framework of G\"odel's universe, the behaviour, and mainly the asymptotic behaviour, of the so-called relativistic diffusion. \par 
    The relativistic diffusion was introduced by J. Franchi and Y. Le Jan in [F-LJ], in the framework of general relativity, on an arbitrary Lorentz manifold, as the only diffusion which is invariant under Lorentz isometries. In this sense, it is the Lorentzian analogue of the Brownian motion on a Riemannian manifold. It lives in fact on the pseudo-unit tangent bundle of the considered Lorentz manifold, and is roughly the integral of Brownian motion of the unit pseudo-sphere of the tangent space. It can also be seen as a random perturbation of the timelike geodesic flow.   

   This article begins with a detailed study of timelike and lightlike geodesics, in a different and more complete way as Kundt [K] did. As a conclusion of the study of lightlike geodesics, a definition (Definition \ref{def.convlum}) of light ray (or boundary point, as an equivalence class of lightlike geodesics, without use of causality) and of convergence to a light ray is given, which appears to be rather natural in this non-causal universe (it can be reinforced to a certain extent : see Remark \ref{rem.conv}). Thus the set of light rays has a natural structure of 3-dimensional boundary, on which the isometry group of G\"odel's universe does operate.  \par 
   Then the relativistic diffusion of G\"odel's universe is introduced. In order to study such a 7-dimensional diffusion, some sub-diffusions are considered, of dimensions 1, 2, and 4.  A leading concern is here to bring out all asymptotic variables of the relativistic diffusion, or in other words, the tail $\sigma$-field of its natural filtration (which is in turn closely related to the Poisson boundary of the relativistic diffusion). The clue in this direction is the idea, suggested and partially worked out in [F-LJ] in the case of Schwarzschild solution, and very recently established in [B] in the flat case of Minkowski space, that convergence to a light ray should eventually occur. The following statement, progressively established below, shows that this general guess stands out as reinforced, also in the case of a non-empty (and even non causal) universe.  \par 
\medskip 
The main results of the present article are summarised in the following. \par \medskip 
\noindent \ub{\bf THEOREM}  \label{the.convlum} \   $(i)$   The relativistic diffusion  is irreducible (on its 7-dimensional phase-space). 
\par 
   $(ii)$ \  Almost surely, the relativistic diffusion path possesses a 3-dimensional asymptotic random variable, and converges to a light ray (in the sense of Definition \ref{def.convlum} and Remark \ref{rem.conv}). \par
   $(iii)$ \ The support of possible light rays the relativistic diffusion can converge to, is the whole 3-dimensional boundary space of light rays.  \par 
\medskip 

   As a consequence of this theorem, and on the basis of some secondary results and considerations, the following conjecture appears to hold likely : by the showing up of the 3-dimensional asymptotic random variable evoked in $(ii)$ above, the whole tail $\sigma$-field of the relativistic diffusion of G\"odel's universe, and then its whole Poisson boundary,  has been brought out. \par 
   The only relativistic case in which the analogous statement has been proved, up to now, is Minkowski space, in [B], by two different methods : Doob's $h$-process conditioning and then couplings, making use of an explicit expression of the laws of already found asymptotic variables ;  or alternatively :  study of the random walk associated with a lifted relativistic diffusion on some (PoincarŽ) fixed locally compact group. The use of both methods does not seem to be easy in the present curved case (likely as in any other curved case), since neither are explicit the laws of the asymptotic variables, nor appears any PoincarŽ-like group. 
   
\eject 

\section{G\"odel's  pseudo-metric} \label{sec.S} 

\begin{defi}\label{def.metr} \   The G\"odel's universe is the manifold $\,\R^4$, endowed with  
coordinates $\,\xi:=(t, x, y, z)$, and with the pseudo-metric (having signature $(+,-,-,-)$) defined by : 
$$ ds^2\, :=\, dt^2 - dx^2 +\,\5\, e^{2\sqrt{2}\,\omega\,x}\,dy^2 +\,2\, e^{\sqrt{2}\,\omega\,x}\,dt\,dy - dz^2\,   , $$ 
for some strictly positive constant $\,\omega\,$. 
\end{defi} 

   The inverse matrix of this pseudo-metric $(\!(g_{ij})\!)$ is as follows : 
$$ (\!(g^{ij})\!) = \pmatrix{ -1 & 0 & 2\, e^{-\sqrt{2}\,\omega\,x} & 0 \cr 0 & -1 & 0 & 0 \cr 
2\, e^{-\sqrt{2}\,\omega\,x} & 0 & -2\, e^{-2\sqrt{2}\,\omega\,x} & 0 \cr 0 & 0 & 0 & -1 \cr } . $$ 

   The unit pseudo-norm relation, defining proper time $\,s\,$, is : 
$$ (\theform)\qquad 1 + \dot t_s^2 +  \dot x_s^2 + \dot z_s^2\, = \, \,\5\,\Big[ e^{\sqrt{2}\,\omega\,x_s}\, \dot y_s + 2\, \dot t_s\Big]^2 . $$ 

   The isometry group of G\"odel's universe is the five-dimensional Lie group generated by : \parn 
1) \  the translations $(t,x, y, z)\mapsto (t+t_0, x, y+y_0, z+z_0)$ of the linear $(t, y, z)$ $3$-subspace  ; \parn
2) \  the hyperbolic dilatations $(t,x, y, z)\mapsto (t, x+x_0, y\,e^{-\sqrt{2}\,\omega\,x_0}, z)$ ;\parn
3) \  the rotational symmetries $\,(u,r, \phi, z)\mapsto (u,r, \phi+\phi_0, z)$, in the new coordinates system   $(u,r, \phi, z)\in \R\times\R_+\times(\R/{2\pi\over\omega}\Z)\times\R\,$ defined by $\, |t-u|< \pi/\omega\,$ and : 
$$ e^{\sqrt{2}\,\omega\,x} = \ch\!(2r) + \sh\!(2r) \cos(\omega\phi) ;\, 
e^{\sqrt{2}\,\omega\,x} \omega\,y = \sh\!(2r) \sin(\omega\phi) ;\, \tg[{\ss{\omega\over 2}}(\phi+t-u)] = e^{-2r} \tg[{\ts{\omega\,\phi\over 2}}] ; $$ 
\centerline{we have indeed \quad $ ds^2 = [du+2\,\sh^2r\, d\phi]^2-2\omega\2 dr^2 -\5\,\sh^2(2r) d\phi^2  - dz^2\,$.} \par 

  G\"odel ([G1], section 4) proved that these three types of isometries generate indeed the full isometry group. As the action of this group is clearly transitive on $\R^4$, G\"odel's universe is an homogeneous space-time.  \par 
  Letting $\,\omega\,$ go to $0$, we recover Minkowski space as limit of G\"odel's universe. 

\subsection{Timelike geodesics} \label{sec.geod} \indf 
   Geodesics are associated with the Lagrangian $\,L(\dot\xi , \xi)$, given by :   
$$ 2\,L(\dot\xi_s,\xi_s) = \dot t_s^2 -  \dot x_s^2 + \,\5\, e^{2\sqrt{2}\,\omega\,x_s}\, \dot y_s^2 + 2\, e^{\sqrt{2}\,\omega\,x_s}\,\dot t_s\,\dot y_s - \dot z_s^2\, .  $$ 

   The equation of geodesics $\quad {\ds {\p\over \p s} \Big(\, {\p L(\dot\xi_s , \xi_s) \over \p \dot \xi^j_s}\Big) = {\p L(\dot\xi_s , \xi_s) \over \p \xi^j_s} }\quad$ reads here : 
$$ (1)\quad \dot t_s + e^{\sqrt{2}\,\omega\,x_s}\, \dot y_s \, = \,  a \; ; \quad (2)\quad e^{2\sqrt{2}\,\omega\,x_s}\, \dot y_s + 2\, e^{\sqrt{2}\,\omega\,x_s}\,\dot t_s \, = \, b  \; ; \quad 
(3)\quad \dot z_s \, =\, c\; ; $$ 
$$ (4)\qquad \ddot x_s + (\omega/\sqrt{2}\,)\, e^{2\sqrt{2}\,\omega\,x_s}\, \dot y_s^2 + \sqrt{2}\,\omega\, e^{\sqrt{2}\,\omega\,x_s}\,\dot t_s\,\dot y_s \, =\, 0\, ; $$ 
for constant \ $\, a,\, b,\, c\,$. \par \medskip 
   
   Equations $(1)$ and $(2)$ jointly are equivalent to : 
$$ (1')\qquad \dot t_s = b\, e^{-\sqrt{2}\,\omega\,x_s} - a\; ; \qquad (2')\qquad \dot y_s = 2\,a\, e^{-\sqrt{2}\,\omega\,x_s} - b\, e^{-2\sqrt{2}\,\omega\,x_s} \; ; $$    
and then using Equations $(1'), (2'), (3)$, we see that Equation $(0)$ is equivalent to : 
$$ 1 + \Big[ b\, e^{-\sqrt{2}\,\omega\,x_s} - a\Big]^2 +  \dot x_s^2 +  c^2\, = \, \,\5\,b^2\, e^{-2\sqrt{2}\,\omega\,x_s} \, ,  $$ 
or equivalently : 
$$ (0')\qquad \dot x_s^2 + \5\, \Big[ 2\, a - b\, e^{-\sqrt{2}\,\omega\,x_s}\Big]^2 =\, a^2 - c^2 -1\, .  $$ 
\par \smallskip  
Note that necessarily \  $\,a^2 \ge 1+ c^2\,$, \  and \  $\,a\,b > 0\,$.   \par\medskip 

   Then, owing to  Equation $(2)$, Equation $(4)$ is equivalent to : 
$$ (4') \qquad {\sqrt{2} \over \omega\, b}\,\,\dot x + y = Y\, , \quad \hbox{ for some constant } \, Y . $$

   Note that Equations $(1)$ and $(2)$ imply : 
$$ (1'')\quad \ddot t_s + \sqrt{2}\,\omega\,e^{\sqrt{2}\,\omega\,x_s}\, \dot x_s\,\dot y_s + 2\, \sqrt{2}\,\omega\, \dot t_s\,\dot x_s \, =\, 0\; ;\; \quad (2'')\quad \ddot y_s - 2\, \sqrt{2}\,\omega\,e^{-\sqrt{2}\,\omega\,x_s}\, \dot t_s\, \dot x_s \, =\, 0\; , $$ 
so that, using the derivative of unit pseudo-norm Relation $(0)$, we see that Equation $(4)$ is implied by Equations $(0), (1), (2), (3)$, unless $\,\dot x_s\equiv 0\,$.  \par\medskip 
    
   Setting  \  $R := a \sqrt{1- (1+c^2)/a^2}$ \   and \   $ k:= {R/ (\sqrt{2}\, a)}\in [0, {\ts{1\over \sqrt{2}}}]$,  \   we must have by $(0')$ : 
$$ \dot x_s = R\,\cos(\omega\, \f_s) \; ,\quad b\, e^{-\sqrt{2}\,\omega\,x_s} = 2\, a - \sqrt{2}\,R\,\sin(\omega\, \f_s) \, , $$    
for some angular component $\,\f_s\,$, whence : 
$$ \dot \f_s\, = \, b\, e^{-\sqrt{2}\,\omega\,x_s}\, =\, 2\,a - \sqrt{2}\,R\,\sin(\omega\, \f_s)  ,  $$ 
and then : 
$$  2\,a\,\omega\, (s- s_0) = \int^{\f_s} {\omega\,d\f \over 1- k\,\sin(\omega\,\f)} = 
{2\over \sqrt{1-k^2}} \, {\rm Arctg}\Big[ {\tg(\omega\,\f_s/2) - k \over \sqrt{1-k^2}} \Big] \, . $$ 
Therefore 
$$ \tg(\omega\,\f_s/2) =  \sqrt{1-k^2}\,\, \tg\Big[ a\,\sqrt{1-k^2}\,\omega\, (s- s_0)  \Big] + k \, , $$     
and 
$$ (5)\qquad  e^{-\sqrt{2}\,\omega\,x_s}\, =\,  { 2\,a\over b}\times \Bigg( 1 - { 2\,k\, \Big( \sqrt{1-k^2}\,\, \tg\Big[ a\,\sqrt{1-k^2}\,\omega\, (s- s_0)  \Big] + k\Big)  \over 1+ \Big( \sqrt{1-k^2}\,\, \tg\Big[ a\,\sqrt{1-k^2}\,\omega\, (s- s_0)  \Big] + k\Big)^2} \Bigg) , $$ 
or equivalently : 
$$ (5')\qquad x_s\, =\,   { -1\over \sqrt{2}\,\omega}\, \log \Bigg[ { 2\,a\over b}\times  
{ (1-k^2)\,\Big( 1+ \tg^2\Big[ a\,\sqrt{1-k^2}\,\omega\, (s- s_0)  \Big]\Big) \over 1+ \Big( \sqrt{1-k^2}\,\, \tg\Big[ a\,\sqrt{1-k^2}\,\omega\, (s- s_0)  \Big] + k\Big)^2} \Bigg] . $$ 

Moreover by $(1')$ we have : 
$$  \dot t_s\, =\, a - \sqrt{2}\,R\,\sin(\omega\, \f_s) = \, \dot\f_s -\,a \, , $$ 
whence 
$$ (6)\qquad  t_s \,=\,  T_0 - a(s-s_0) + {{2\over \omega}}\, {\rm Arctg}\Big( \sqrt{1-k^2}\,\, \tg\Big[ a\,\sqrt{1-k^2}\,\omega\, (s- s_0)  \Big] + k \Big) . $$ 
In this formula $(6)$, the successive determinations of $\,{\rm Arctg}\,$, at the successive values \  $s\in s_0+ {\pi\over a\,\sqrt{1-k^2}\,\omega}\,\Z\,$, are understood to be chosen conveniently, in order that the absolute time coordinate $(t_s)$ be continuous, as it must be. \   
Observe that $(t_s)$ is strictly monotonic if and only if  $\,k\le \5\,$, or equivalently, if and only if  $\,(1+c^2)\le a^2\le 2(1+c^2)$. \par \smallskip  

   Finally,  by $(4')$ and $(5')$ we have : 
$$ (7)\qquad  y_s\, = \,  Y + { 2\,a\,k\over b\,\omega} - { 4\,a\,k/ (b\,\omega) \over 1+ \Big( \sqrt{1-k^2}\,\, \tg\Big[ a\,\sqrt{1-k^2}\,\omega\, (s- s_0)  \Big] + k\Big)^2}\, , $$ 
which is consistent with  $(2')$ : \  
${\ds  \dot y_s\, =\,4\,{a^2k\over b}\, [1-k\,\sin(\omega\, \f_s)] \sin(\omega\, \f_s)}$. \par 
\if{
$$ =  8\,{a^2k\over b} ({1-k^2})\, { \bigg[\sqrt{1-k^2}\, \tg\Big[ a\sqrt{1-k^2}\,\omega (s- s_0) \Big] + k \bigg] \! \times\! \bigg[ 1+ \tg^2\Big[ a\sqrt{1-k^2}\,\omega (s- s_0) \Big]\bigg]   \over \bigg[1+ \Big( \sqrt{1-k^2}\,\, \tg\Big[ a\,\sqrt{1-k^2}\,\omega\, (s- s_0)  \Big] + k\Big)^2\bigg]^2} \, , $$ 
}\fi 
\par\medskip 

   Observe from Equations $(5)$ and $(7)$ that every timelike geodesic has a bounded, periodic projection in the $(x,y)$-plane, and moreover that it obeys the following relation : 
$$  (8)\qquad  \Big[ {b\over 2\, a}\, e^{-\sqrt{2}\,\omega\,x_s} - 1\Big]^2 + \Big[{\omega\, b\over 2\, a}\, (y_s - Y) \Big]^2 =\, k^2\, . $$ 

\brem \label{rem.ptpart} \   {\rm  The case $\,k=0\,$ is particular. It implies (using Equations $(5), (0'), (1'), (2')$) : \  $\dot t^2 = 1+ \dot z^2\,$ and $\,\dot x = \dot y = 0\,$, \  and then : \par\smallskip 
\centerline{ $(x_s,y_s)$ constant and $\,t_s = t_0+ a\,s\,$, $\,z_s = z_0+ c\,s\,$, with $\, a^2= 1+c^2\,$.} \parn 
Reciprocally, if \  $\dot x_0 = \dot y_0 = 0\,$, then (by Equations $(2'), (0')$) the corresponding geodesic must satisfy also $\,k =0\,$, and then be included in the phase subspace defined by :  
$$ \EE_0 := \Big\{\dot x = \dot y =0\Big\} = \Big\{ \dot t^2 = 1+ \dot z^2\,\,;\; \dot x = \dot y =0\Big\} . $$ 
Therefore, the case $\,k=0\,$ corresponds to the geodesically stable phase subspace $ \EE_0\,$. 
}\erem 

\brem \label{rem.causal} \quad   {\rm  Every timelike geodesic is defined for all proper times $\,s\,$, unbounded and causal. Moreover, it never accumulates near its past. Indeed, if for proper times $\,s<s'\,$ we had $\,x_{s'}=x_s\,$, then by Equation $(5)$ we should have : 
$$ \tg\Big[ a\sqrt{1-k^2}\,\omega (s'- s_0)\Big]=\tg\Big[ a\sqrt{1-k^2}\,\omega (s- s_0)\Big] , \; \hbox{ then } \; s' = s+ {\ts{\pi\, n\over |a|\omega\sqrt{1-k^2}}}\; \hbox{ with } \; n\in\N^*, $$ 
whence by Equation $(6)$ : \  $\,t_{s'}-t_s = {\rm sign}(a)\Big( {{2\pi\, n\over \omega}} - {\pi\, n\over \omega\sqrt{1-k^2}}\Big)$, \  and then \  
$$ |t_{s'}-t_s| = {n\pi\over \omega}\,(2-(1-k^2)^{-1/2}) \,\ge\, {\pi\over \omega}\, . $$ 
}\erem

  The following statement, which will be used later to ensure the irreducibility of the relativistic diffusion, shows up the non-causal structure of G\" odel's universe, despite the preceding remark \ref{rem.causal} : the causal past of any point  of G\" odel's universe is the whole G\" odel's universe. In particular, the causal boundary, in the sense of Penrose (or Geroch-Kronheimer-Penrose, see ([H-E], section 6.8)), reduces to a single point. 
\bpro \label{pro.geod} \  The G\" odel's universe $\Big( \R^4,((g_{ij}))\Big)$ is geodesically transitive : any two points of it can be linked by a piece-wise geodesic timelike continuous path. 
\epro 
\ub{Proof} \quad Observe from Remark \ref{rem.ptpart} that (taking $\,k=c=0$) there are timelike geodesics moving at will the coordinate $\,t\,$, without changing any other coordinate, and  that (taking $\,k=0,\,c\not=0 $) there are timelike geodesics moving at will the coordinate $\,z\,$, without changing the coordinates $(x,y)$. \par 
   Hence it will be sufficient to move piece-wise geodesically the coordinates $(x,y)$, forgetting henceforth the coordinates $(t,z)$.   \par 
   Observe that any given geodesic can move the coordinate $\,x\,$ only by a uniformly bounded value, since by Equation $(5)$ we have : 
$$ e^{-\sqrt{2}\,\omega\,(x_s-x_{s'})}\, =\, { [1+ \tg^2(\omega\, \f_s)] \Big(1+[\sqrt{1-k^2}\,\, \tg(\omega\, \f_{s'}) + k]^2\Big) \over  \Big(1+[\sqrt{1-k^2}\,\, \tg(\omega\, \f_{s}) + k]^2\Big) [1+ \tg^2(\omega\, \f_{s'})]} \, , $$ 
so that several geodesic arcs are needed. \   Fix coordinates $(x,y)$, which we want to move to other fixed coordinates say $(x',y')$. By a finite number of geodesic moves, according to the equation displayed just above, we can move $(x,y)$ to $(x',y'')$, for some $\,y''$. \par
   Then, the quotient in the equation displayed just above equals 1 as soon as   $$\f_s+\f_{s'}\, =\, \omega\1\,\arctg\Big[(2k\sqrt{1-k^2}\,)\Big/(1+k^2-\sqrt{1-k^2}\,)\Big] \, , $$
so that we can choose a geodesic arc, independently of the value of the parameter $\,b\,$,  such that this holds and such that \  \mbox{$[\sqrt{1-k^2}\,\, \tg(\omega\, \f_{s'}) + k]^2\not= [\sqrt{1-k^2}\,\, \tg(\omega\, \f_{s}) + k]^2$}. \  Finally, on such geodesic, by Equation $(7)$ we have : 
$$ y_s-y_{s'}\, =\, {4\,a\,k\over b\,\omega}\times\bigg[ { 1 \over  1+[\sqrt{1-k^2}\,\, \tg(\omega\, \f_{s'}) + k]^2} - { 1 \over  1+[\sqrt{1-k^2}\,\, \tg(\omega\, \f_{s}) + k]^2}\bigg]  , $$ 
proving that, choosing conveniently the parameter $\,b\,$, we can move $(x',y'')$ to $(x',y')$. $\;\diamond$ 

\if{
\subsection{Geodesics from a fixed observer} \label{sec.fixob} \indf 
   Viewed by a fixed observer, geodesics are parametrised by absolute time $\,t\,$, and then have trajectories \  $(X_t,Y_t,Z_t) := (x_{S(t)} , y_{S(t)}, z_{S(t)})$, where the function $\,S\,$ is determined by the relation \  $ s \, = \, S(t_s)$. By Equation $(6)$, we have : 
   }\fi 
 
\subsection{Lightlike geodesics} \label{sec.ngeod} \indf 
   Equations $(1), (2), (3), (4)$ remain the same, while the pseudo-norm equation $(0)$ is replaced by : 
$$ (0'')\qquad \dot t_s^2 +  \dot x_s^2 + \dot z_s^2\, = \, \,\5\,\Big[ e^{\sqrt{2}\,\omega\,x_s}\, \dot y_s + 2\, \dot t_s\Big]^2 . $$ 
As previously, knowing Equations $(1)$ and $(2)$, Equation $(4)$ is equivalent to
$$ (4') \qquad {\sqrt{2} \over \omega\, b}\,\,\dot x + y = Y\, , \quad \hbox{ for some constant } \, Y , $$
and is implied by Equations $(0''), (1), (2), (3)$.  Thus  lightlike geodesics are the solutions to the system :
$$ (1')\quad \dot t_s = b\, e^{-\sqrt{2}\,\omega\,x_s} - a\; ; \quad (2')\quad \dot y_s = 2\,a\, e^{-\sqrt{2}\,\omega\,x_s} - b\, e^{-2\sqrt{2}\,\omega\,x_s} \; ; \quad 
(3)\quad \dot z_s \, =\, c\; ; $$    
and 
$$ (0''')\qquad \dot x_s^2 + \5\, \Big[ 2\, a - b\, e^{-\sqrt{2}\,\omega\,x_s}\Big]^2 =\,  a^2 - c^2 \, .  $$ 
Hence we must have again \  $ a^2 \ge c^2\,$, and $\,a b > 0\,$ (if not, the trajectory must be constant), and, setting now \  
$$\kappa\, :=\, \sqrt{\5 (1- c^2/a^2)} \in [0, {\ts{1\over \sqrt{2}}}]\, , $$
we get the same equations $(5), (6), (7), (8)$ as above, merely with $\,\kappa\,$ replacing $\,k\,$. \par \smallskip 

   As the parameter $\,s\,$ cannot any longer have here a meaning like proper time, but stands only for an affine parameter, determined up to a change $\,s\mapsto u\, s+v\,$, the constant $(a,b,c)$ is now irrelevant. Note that on the contrary, the constant $Y$ of Equation $(4')$ is consistent. The only meaningful a priori parameter for a lightlike geodesic is the  ``impact'' parameter : \  
$$ (9) \qquad B = (\ell,\rr,Y) := \Big({c\over a},{b\over a}, Y \Big) \in \BB := [-1,1]\times \R^*_+\times \R\, . $$ 
   Eliminating $\,s\,$, we see indeed that a lightlike geodesic having impact parameter $\,B\,$ solves : \par
\vbox{ 
$$ \ell\, dt\, =\, (\rr\, e^{-\sqrt{2}\,\omega\,x} - 1)\, dz \;;\quad (2-\rr\, e^{-\sqrt{2}\,\omega\,x})  \, e^{-\sqrt{2}\,\omega\,x}\, dt \,=\, (\rr\, e^{-\sqrt{2}\,\omega\,x} - 1)\, dy \;;\quad  $$ 
$$ (\rr\, e^{-\sqrt{2}\,\omega\,x} - 1)\, dx\, =\,\pm \sqrt{ 1-\ell^2-\5 (2-\rr\, e^{-\sqrt{2}\,\omega\,x}) ^2} \,\, dt\, . $$ } \par \medskip

   Let us sum up the description of lightlike geodesics in the following statement. 
\bpro \label{pro.geodlum} \quad  Any lightlike geodesic $(x_\tau, y_\tau, z_\tau, t_\tau)$ having impact parameter \parn $B= (\ell,\rr,Y)\in\BB\,$ satisfies, for an additional parameter $(Z_0,T_0)\in \R^2$ and for any real $\,\tau$ : 
$$ e^{-\sqrt{2}\,\omega\,x_\tau}\, =\,  { 2\over \rr}\times \Bigg( 1 - { 2\,\sqrt{1-\ell^2}\, \Big( \sqrt{1+\ell^2}\,\, \tg\tau + \sqrt{1-\ell^2}\,\Big)  \over 2+ \Big( \sqrt{1+\ell^2}\,\, \tg\tau + \sqrt{1-\ell^2}\Big)^2} \Bigg) ; $$ 
$$  y_\tau\, = \, {Y} + {\sqrt{2(1-\ell^2)}\over \omega \rr} \Bigg( 1 - { 4 \over 2+ \Big( \sqrt{1+\ell^2}\,\, \tg\tau + \sqrt{1-\ell^2}\, \Big)^2}\Bigg) \, ; $$ 
$$ z_\tau\, =\, Z_0 + {\ell\, \tau\over \omega\sqrt{(1+\ell^2)/2}} \, ; $$ 
$$   t_\tau \,=\,  T_0 - {\tau\over \omega\sqrt{(1+\ell^2)/2}} + {{2\over \omega}}\, {\rm Arctg}\bigg( \sqrt{(1+\ell^2)/2}\,\, \tg\tau + \sqrt{(1-\ell^2)/2}\, \bigg) ; $$ 
$$ \CC_{B} : \quad  \Big[ {\rr\over 2}\, e^{-\sqrt{2}\,\omega\,x_\tau} - 1\Big]^2 + \Big[{\omega\, \rr\over 2}\, (y_\tau - Y)\Big]^2 =\, {1-\ell^2\over 2}\,  . $$ 
\epro

\brem \label{rem.geodlum} \  {\rm  The last equation in Proposition \ref{pro.geodlum} shows that to any given lightlike geodesic is associated a cylinder $\CC_{B}$, parallel to the $(t,z)$-coordinate plane.  Reciprocally, by Proposition \ref{pro.geodlum} again, any lightlike geodesic which is drawn on the cylinder $\CC_{B}$ has a prescribed projection on the $(x,y)$-coordinate plane (up to  changing affine parameter $\tau$). The equations displayed in Proposition \ref{pro.geodlum} define a lightlike geodesic associated to any given $\,B\in\BB$.\par 
   Considering then any  continuous angular parameter $\,\f = \f_\tau$ (determined modulo $2\pi/\omega$) such that :
$$ \tg (\omega\,\f_\tau/2)\,  = \, \sqrt{(1+\ell^2)/2}\;\, \tg\tau + \sqrt{(1-\ell^2)/2}\,\, , $$ 
then by Proposition \ref{pro.geodlum} we have :
$$  \rr\, e^{-\sqrt{2}\,\omega\,x_\tau}\, =\, 2 - \sqrt{2\,(1-\ell^2)}\,\sin (\omega\,\f_\tau) \quad \hbox{ and } \quad \omega\, \rr\; y_\tau \,=\, \omega\, \rr\,Y - \sqrt{2\,(1-\ell^2)}\,\cos (\omega\,\f_\tau) , $$ 
together with : 
$$  t_\tau \,=\,  T_0 - {\tau\over \omega {\ss\sqrt{(1+\ell^2)/2}}} + \f_\tau \; ,\quad z_\tau + \ell\, t_\tau\, =\, Z_0 + \ell\,T_0 +  \ell\, \f_\tau\, . $$ 

   Since the function $\,\tau\mapsto \omega\, \f_\tau - 2\, \tau\,$ is $\,\pi$-periodic, the functions  \  $$\tau\mapsto  \omega\, t_\tau - 2\, (1-[2(1+\ell^2)]^{-1/2})\, \tau \quad \hbox{and} \quad  \tau\mapsto z_\tau - {\ell\, t_\tau\over \sqrt{2(1+\ell^2)}\,-1} = {\ts Z_0 - {\ell ( T_0+ \f_\tau- 2\tau/\omega) \over \sqrt{2(1+\ell^2)}\,-1} } $$ 
are $\,\pi$-periodic too. This implies that $\,t_\tau\,$ wanders out to infinity, nearly linearly, and that the projection on the $(t,z)$-coordinate plane of any lightlike geodesic has an asymptotic direction : 
$$ \lim_{\tau\to\pm\infty}\, { z_\tau\over t_\tau}\, =\,  {\ell\over \sqrt{2(1+\ell^2)}\,-1}\, . $$ 

    This prescribes geometrically the sign of parameter $\,\ell\,$, which was not determined by the cylinder $\CC_{B}$ alone, which however determines $(|\ell|,\rr, Y)$. Note that the impact parameter $ B = (\ell,\rr, Y)\in \BB\,$ has thus indeed a clear geometrical meaning.  \par 
    Note finally that the additional parameter $(Z_0,T_0)$ depends on a translation on the parameter $(\tau,\f_\tau)$, and then, contrary to $\,B=(\ell,\rr, Y)$,  is geometrically irrelevant. }\erem 
 \medskip   
   
   Recall that in a strongly causal space-time, it seems natural to use the causal boundary, in the sense of Penrose, to classify lightlike geodesics by gathering in an equivalence class, called a light ray, all geodesics which converge to a given causal boundary point (having asymptotically the same past, see ([H-E], section 6.8)).  On the contrary, in the present setting (recall Proposition \ref{pro.geod}), such classification is totally inoperative. It seems that no alternative classification has been proposed so far, which be relevant in a non-causal setting.  \par
   However, owing to the above remark \ref{rem.geodlum}, we are led to adopt here the following alternative classification of lightlike geodesics into light rays, and then also, to see the 3-dimensional space of light rays as an alternative notion of (non-causal) boundary, as follows. 

\bdefi \label{def.convlum} \  Let us call \ub{light ray}, or \ub{boundary point}, of G\"odel's universe,  any equivalence class of lightlike geodesics, identifying those which have the same impact  parameter $\,B=(\ell,\rr, Y)\in \BB\,$. \   Thus $\,\BB = [-1,1]\times \R^*_+\times \R\,$ is the  \ub{boundary} of G\"odel's universe. \par 
   Let us say that a path $\,s\mapsto \xi_s =  (t_s,x_s,y_s,z_s)$ of class $\,C^1$ in G\"odel's universe \parn 
\ub{converges to the light ray} $\,B=(\ell,\rr, Y)$ if, setting :   
$$  a_s\, :=\, \dot t_s + e^{\sqrt{2}\,\omega\,x_s}\, \dot y_s \; , \quad  b_s\, :=\, e^{\sqrt{2}\,\omega\,x_s}\,(2\, \dot t_s + e^{\sqrt{2}\,\omega\,x_s}\, \dot y_s) \, , \;  \hbox{ and } \quad  Y_s\, :=\, {\sqrt{2}\,\,\dot x_s \over \omega\, b_s} + y_s\; , $$ 
the following convergences hold, as $\,s\to +\ii$ : 
$$ {\dot z_s\over a_s} \lra \ell \;,\quad {b_s\over a_s} \lra \rr\; ,\quad Y_s  \lra Y\, ,\quad 
 \Big[ {\rr\over 2}\, e^{-\sqrt{2}\,\omega\,x_s} - 1\Big]^2 + \Big[{\omega\, \rr\over 2}\, (y_s - Y)\Big]^2 \lra {1-\ell^2\over 2}\, . $$ 
\edefi 
  
     This notion of convergence to the boundary $\,\BB\,$ can be reinforced to a certain extent : see Remark \ref{rem.conv}, concluding this article. \   We saw above that any lightlike geodesic belonging to a light ray $B$ converges to it. On the contrary, a timelike geodesic does not converge to any light ray : by  Section \ref{sec.geod}, we get indeed $\,B=({c\over a},{b\over a},Y)\in\BB\,$, but the cylinder $\CC_B\,$ has to small a ``radius", since by Equation $(8)$ we have $\,k^2 = [1-\ell^2 - a\2]/2  <  [1-\ell^2]/2\,$.

\bpro \label{pro.opB} \  The isometry group of G\"odel's universe (recall Section \ref{sec.S}) does operate on the boundary $\,\BB$, so that  the above definition \ref{def.convlum}  is consistent. \  Precisely, it works as follows. \parn  
1) \  The translation $(t,x, y, z)\mapsto (t+t_0, x, y+y_0, z+z_0)$ changes $(\ell,\rr, Y)$ into $(\ell,\rr, Y+y_0)$. \parn
2) \  The hyperbolic dilatation $(t,x, y, z)\mapsto (t, x+x_0, y\,e^{-\sqrt{2}\,\omega\,x_0}, z)$  changes $(\ell,\rr, Y)$ into $$(\ell,\rr\, e^{\sqrt{2}\,\omega\,x_0}, Y e^{-\sqrt{2}\,\omega\,x_0})\, . $$ 
3) \  The rotational symmetry $\,(u,r, \phi, z)\mapsto (u,r, \phi+\phi_0, z)$ changes $(\ell,\rr, Y)$ into $$ \left( \ell\, , \, {\alpha} + [\rr - {\alpha}]\cos(\omega\,\phi_0) - \omega\,\rr \, Y\, \sin(\omega\,\phi_0)\, ,\, { \omega\,\rr\,\cos(\omega\,\phi_0)+ [\rr -{\alpha}] \sin(\omega\,\phi_0) \over  \omega\,[{\alpha}+ [\rr -{\alpha}]\cos(\omega\,\phi_0)- \omega\,\rr\,\sin(\omega\,\phi_0)]  }\right) ,  $$
where \  $\alpha := {2a\,\ch(2r) -\sh^2(2r)\,\dot\phi\over \dot u + 2\,\sh^2r\,\,\dot\phi}\, $ is constant under $\,\phi\mapsto\phi+\phi_0\,$, and on each geodesic. We have indeed : \   
$\alpha = {\rr\over 2}\, (1+\omega^2Y^2)+  {1+\ell^2\over \rr}\,$, or equivalently :  \    
$\rr-\alpha = {\rr\over 2}\, (1-\omega^2Y^2) -  {1+\ell^2\over \rr}\,$. 
\epro 
\ub{Proof} \quad  The two first items are straightforward. On the other hand, the action of the rotational symmetry $\,(u,r, \phi, z)\mapsto (u,r, \phi+\phi_0, z)$ is not so obvious. However, a computation shows that we have in coordinates $(u,r, \phi, z)$ : \  $\, a = \dot u + 2(\sh r)^2\,\dot\phi\,$, \  and : 
$$ b = A + \Psi\, \cos(\omega\,\phi) -2\omega\1\dot r\,\sin(\omega\,\phi) \, , \quad 
Z:= \omega\,b\, Y = \Psi\, \sin(\omega\,\phi) + 2\omega\1\dot r\,\cos(\omega\,\phi)\, ,  $$ 
with 
$$ A := 2a\,\ch(2r) -\sh^2(2r)\,\dot\phi \quad \hbox{ and } \quad \Psi := [2a -\ch(2r)\,\dot\phi]\,\sh(2r)\, . $$ \parn 
Note that \  $A = 2a+4[a -\ch^2r\,\dot\phi]\,\sh^2r = 2a+2\,{\partial\over \partial\dot\phi} L$ \   is seen to be constant on each geodesic, by looking at the expression of the Lagrangian $\,L\,$ in coordinates $(u,r,\phi,z)$. Alternatively, a computation yields : \  
$ 2A = b + \omega^2 b\,(2 Y y - y^2) + (4a - b\,e^{-\sqrt{2}\,\omega\,x})\, e^{-\sqrt{2}\,\omega\,x}$, whence by using Proposition \ref{pro.geodlum} : 
$$ 2\alpha -\rr =  {\ts{2A\over a}} -\rr = \omega^2\rr\! \left[ Y^2 - {\ts{2(1-\ell)^2\over \omega^2\rr^2}} \! \left[ 1- {4\over 2+\tg^2({\omega\f_\tau\over 2})}\right]^2 \right] + {{1\over\rr}}\! \left[ 4 - \left[  {4\sqrt{1-\ell^2}\,\tg^2({\omega\f_\tau\over 2}) \over 2+\tg^2({\omega\f_\tau\over 2})}\right]^2 \right]
$$ 
$$ =\, \omega^2 \rr\, Y^2 + {4\over\rr} - 2\, {(1-\ell^2)\over\rr} \; , \hbox{ whence } \quad 
\alpha = {\rr\over 2}\, (1+\omega^2Y^2)+ {1+\ell^2\over \rr}\, . $$ 
 
   Now we have at once : under $\,\phi\mapsto\phi+\phi_0\,$, \  $(a,b,Z)$ is changed into 
$$ \Big( a\, ,\,  A+ [b-A]\cos(\omega\,\phi_0)-Z\sin(\omega\,\phi_0)\, ,\, Z\cos(\omega\,\phi_0)+ [b-A]\sin(\omega\,\phi_0)\Big) , $$
so that $(\ell,\rr, Y)$ is changed into (recall from the above that $\,\alpha = A/a$) : 
$$ \left( \ell\, , \, {\ts{A\over a}} + [\rr - {\ts{A\over a}}]\cos(\omega\,\phi_0) - \omega\,\rr \, Y\, \sin(\omega\,\phi_0)\, ,\, { \omega\,\rr\,\cos(\omega\,\phi_0)+ [\rr -{\ts{A\over a}}] \sin(\omega\,\phi_0) \over  \omega\,[{\ts{A\over a}}+ [\rr -{\ts{A\over a}}]\cos(\omega\,\phi_0)- \omega\,\rr\,\sin(\omega\,\phi_0)]  }\right)\! .\;\; \diamond  $$

\subsection{Ricci curvature and energy tensor} \indf 
      Recall that the Christoffel symbols are computed by : \par 
$$ \G_{ij}^k = \5 \, g^{k\ell}\Big( {\p g_{\ell j}\over \p\xi^i} + {\p g_{i \ell}\over \p\xi^j} - {\p g_{i j}\over \p\xi^\ell}\Big)  \, , $$ 
and that the Ricci tensor $\,(R_{ij})\,$ is computed by 
$$ R_{ij} = {\p\G^k_{ij}\over \p\xi^k} - {\p\G^k_{ik}\over \p\xi^j} + \G^k_{k \ell}\G^\ell_{ij} - 
\G^k_{i \ell}\G^\ell_{jk} \, . $$ 

From Equations $(1'), (2'), (3)$, we find all non-vanishing Christoffel coefficients : 
$$ \G_{xy}^t  = \G_{ty}^x = (\omega/\sqrt{2}\,)\, e^{\sqrt{2}\,\omega\,x}\, , \; \G_{tx}^t =  \sqrt{2}\,\omega\, ,  \; \G_{yy}^x = (\omega/\sqrt{2}\,)\, e^{2\sqrt{2}\,\omega\,x}\, , \;  \G_{tx}^y =  - \sqrt{2}\,\omega\,e^{-\sqrt{2}\,\omega\,x}\, . $$  
Therefore we get all non-vanishing Ricci coefficients : 
$$ R_{tt} = 2\omega^2\, ; \; R_{ty} = R_{yt} = 2\omega^2\, e^{\sqrt{2}\,\omega\,x}\, ; \; R_{yy} = 2\omega^2\, e^{2\sqrt{2}\,\omega\,x}\, . $$ 
Hence, the scalar curvature is \quad  $ R\, = \, g^{ij}\,R_{ij}\, = \,2\omega^2 .  $ \par  

   Einstein equations  
$$ R_{ij} - \5\, R\, g_{ij} + \Lambda\,g_{ij}  = \, T_{ij} $$
are satisfied, with cosmological constant $\,\Lambda = \omega^2\,$ representing a positive pressure, and energy tensor $(T_{ij}) = (R_{ij}) = (u_i u_j)$, where \  $u:= (\sqrt{2}\,\omega, 0, \sqrt{2}\,\omega \, e^{\sqrt{2}\,\omega\,x},0)$ represents the four-velocity of matter, which rotates with constant velocity $\,\omega\,$. The energy is thus : 
$$ E(\xi,\dot\xi) :=\,T_{ij}(\xi)\,\dot\xi^i\,\dot\xi^j = 2{\omega^2}\, \Big[ \dot t + e^{\sqrt{2}\,\omega\,x}\, \dot y \Big]^2 = 2{\omega^2}\, a(\xi,\dot\xi)^2\, . $$ 
  
\section{Relativistic diffusion} \label{sec.RD} \indf 
   Recall from [F-LJ] that the general expression of the relativistic operator $\LL\,$ is : 
$$  \LL\, =\, \dot\xi^k\,{\p\over\p\xi^k} + \Big( {\ts{3\,\s^2\over 2}}\, \dot\xi^k - \dot\xi^i\dot\xi^j\, \G_{ij}^k(\xi) \Big) \,{\p\over\p\dot\xi^k} + {\ts{\s^2\over 2}}\, (\dot\xi^k\,\dot\xi^l - g^{kl}(\xi)) \,{\p^2\over\p \dot\xi^k\,\p \dot\xi^l} \, , $$ 
$\sigma\,$ being an arbitrary fixed positive (speed or heat) parameter. \par 

   Equivalently in the present setting, the relativistic diffusion $(\xi_s,\dot\xi_s)$, in coordinates $\,\xi=(t,x,y,z)$, solves the following system of stochastic differential equations : 
$$ dt_s = \dot t_s\, ds \; ; \quad dx_s = \dot x_s\, ds \; ; \quad dy_s = \dot y_s\, ds \; ; \quad dz_s = \dot z_s\, ds \; ; $$     
$$ d \dot t_s\, =\, - 2\, \sqrt{2}\,\omega\, \dot t_s\,\dot x_s\, ds -\sqrt{2}\,\omega\,e^{\sqrt{2}\,\omega\,x_s}\, \dot x_s\,\dot y_s\, ds + {\ts{3\,\sigma^2\over 2}}\, \dot t_s\,ds + \sigma \, dM^t_s \; ; $$ 
$$ d \dot x_s\, =\, - \sqrt{2}\,\omega\, e^{\sqrt{2}\,\omega\,x_s}\,\dot t_s\,\dot y_s \, ds - (\omega/\sqrt{2}\,)\, e^{2\sqrt{2}\,\omega\,x_s}\, \dot y_s^2\, ds + {\ts{3\,\sigma^2\over 2}}\, \dot x_s\,ds + \sigma \, dM^x_s \; ; $$ 
$$ d \dot y_s\, =\, 2\, \sqrt{2}\,\omega\,e^{-\sqrt{2}\,\omega\,x_s}\, \dot t_s\, \dot x_s\, ds + {\ts{3\,\sigma^2\over 2}}\, \dot y_s\,ds + \sigma \, dM^y_s \; ; $$ 
$$ d \dot z_s\, =\, {\ts{3\,\sigma^2\over 2}}\, \dot z_s\,ds + \sigma \, dM^z_s \; ; $$ 
\smallskip \noindent 
where the $\,\R^4$-valued martingale \   $M_s := (M^t_s,M^x_s, M^y_s, M^z_s)$ has (rank 3) quadratic covariation matrix : 
$$ ((K^{ij}_s)) := {\langle dM^i_s, dM^j_s\rangle\over ds} =  \pmatrix{ \dot t_s^2 +1 & \dot t_s\,\dot x_s & \dot t_s\,\dot y_s - 2\, e^{-\sqrt{2}\,\omega\,x_s} & \dot t_s\,\dot z_s \cr 
\dot t_s\,\dot x_s & \dot x_s^2 +1 & \dot x_s\,\dot y_s & \dot x_s\,\dot z_s \cr 
\dot t_s\,\dot y_s - 2\, e^{-\sqrt{2}\,\omega\,x_s} & \dot x_s\,\dot y_s & \dot y_s^2 +2\, e^{-2\sqrt{2}\,\omega\,x_s} & \dot y_s\,\dot z_s \cr \dot t_s\,\dot z_s & \dot x_s\,\dot z_s & \dot y_s\,\dot z_s & \dot z_s^2+1 \cr } . $$ 

    Recall that the unit pseudo-norm relation reads : 
$$ (0)\qquad 1 + \dot t_s^2 +  \dot x_s^2 + \dot z_s^2\, = \, \,\5\,\Big[ e^{\sqrt{2}\,\omega\,x_s}\, \dot y_s + 2\, \dot t_s\Big]^2 . $$ 

   Thus, the relativistic diffusion $(\xi_s,\dot\xi_s)$ is $7$-dimensional, having phase space : 
$$ \EE := \Big\{ (t,x,y,z, \dot t,\dot x,\dot y,\dot z)\in\R^8\,\Big|\,1 + \dot t^2 +  \dot x^2 + \dot z^2\, = \, \,\5\,\Big[ e^{\sqrt{2}\,\omega\,x}\, \dot y + 2\, \dot t\Big]^2 \Big\} ,  $$ 
or equivalently : 
$$ \EE = \Big\{ (t,x,y,z, \dot t,\dot x,\dot y,\dot z)\in\R^8\,\Big|\,1 +  \dot x^2 + \dot z^2 + \5\, e^{2\sqrt{2}\,\omega\,x}\, \dot y^2\, = \, \Big[ \dot t + e^{\sqrt{2}\,\omega\,x}\, \dot y \Big]^2 \Big\} . $$ 

  Note that the particular phase subspace distinguished in Remark \ref{rem.ptpart} :  
$$ \EE_0 := \EE\cap \Big\{\dot x = \dot y =0\Big\} = \EE\cap \Big\{ \dot t^2 = 1+ \dot z^2\,\,;\; \dot x = \dot y =0\Big\} $$  
is clearly not stable under the relativistic diffusion $(\xi_s,\dot\xi_s)$, contrary to the geodesic flow, and even instantly unstable : starting from any point in $\,\EE_0\,$, its exit time from $\,\EE_0\,$ is null. 

\subsection{Reduction of the dimension} \label{sec.RedDim} \indf 
   The study of geodesics induces to consider the following quantities (which, as $\,\dot z_s \,$, are constant  along each geodesic), setting (as in Definition \ref{def.convlum}) : 
$$ (10) \qquad   a_s\, :=\, \dot t_s + e^{\sqrt{2}\,\omega\,x_s}\, \dot y_s \qquad \hbox{ and } \qquad  b_s\, :=\, e^{\sqrt{2}\,\omega\,x_s}\,(2\, \dot t_s + e^{\sqrt{2}\,\omega\,x_s}\, \dot y_s) \, . $$ 

Then we have : 
$$ d a_s\, =\,  {\ts{3\,\sigma^2\over 2}}\, a_s\,ds + \sigma \, dM^a_s =\,  {\ts{3\,\sigma^2\over 2}}\, a_s\,ds + \sigma \, (dM^t_s+ e^{\sqrt{2}\,\omega\,x_s}\,dM^y_s) \; ; $$ 
and 
$$ d b_s\, =\, {\ts{3\,\sigma^2\over 2}}\, b_s\,ds + \sigma \, dM^b_s =\,  {\ts{3\,\sigma^2\over 2}}\, a_s\,ds + \sigma\, e^{\sqrt{2}\,\omega\,x_s}\,(2\, dM^t_s+ e^{\sqrt{2}\,\omega\,x_s}\,dM^y_s) \, . $$ 
Moreover we have : 
$$ d\dot x_s\, =\, (\omega/\sqrt{2}\,)\, e^{-2\sqrt{2}\,\omega\,x_s}\, b_s^2\, ds - \sqrt{2}\,\omega\, e^{-\sqrt{2}\,\omega\,x_s}\,a_s\,b_s\, ds + {\ts{3\,\sigma^2\over 2}}\,\dot x_s\,ds + \sigma\, dM^x_s\, , $$ 
and the $\,\R^4$-valued martingale \   $\tilde M_s := (M^a_s,M^b_s, M^x_s, M^z_s)$ has (rank 3) quadratic covariation matrix : 
$$ ((\tilde K^{ij}_s)) =  \pmatrix{ a_s^2 -1 & a_s\,b_s - 2\,e^{\sqrt{2}\,\omega\,x_s} & a_s\,\dot x_s & a_s\,\dot z_s \cr 
a_s\,b_s - 2\,e^{\sqrt{2}\,\omega\,x_s} & b_s^2 - 2\,e^{2\sqrt{2}\,\omega\,x_s} & b_s\,\dot x_s & b_s\,\dot z_s \cr 
a_s\,\dot x_s & b_s\,\dot x_s & \dot x_s^2 + 1 & \dot x_s\,\dot z_s \cr 
a_s\,\dot z_s & b_s\,\dot z_s & \dot x_s\,\dot z_s & \dot z_s^2+1 \cr } . $$ 
From this, we deduce the following.
\bcor \label{cor.reddiff} \  The ($7$-dimensional) relativistic diffusion $(\xi_s,\dot\xi_s)$ admits the following sub-diffusions : \quad  $(a_s)$ ; \quad  $(\dot z_s)$ ; \quad $(z_s,\dot z_s)$ ; \quad $(a_s,\dot z_s)$ ; \quad $(a_s, z_s,\dot z_s)$ ; \quad $(x_s,\dot x_s,a_s,b_s)$ . 
\ecor 

    The unit pseudo-norm relation can be written :  
$$ (00)\qquad 1 + \dot x_s^2 + \dot z_s^2 + (a_s - e^{-\sqrt{2}\,\omega\,x_s}\, b_s)^2 \, = \, \,\5\,e^{-2\sqrt{2}\,\omega\,x_s}\, b_s^2\, , $$ 
or equivalently : 
$$ (00')\qquad 1 + \dot x_s^2 + \dot z_s^2 + \5\,(2\,a_s - e^{-\sqrt{2}\,\omega\,x_s}\, b_s)^2 \, = \, a_s^2\, . $$ 

   Hence the phase space $\,\EE\,$ of the relativistic diffusion $(\xi_s,\dot\xi_s)$ can be written equivalently : 
$$ \EE = \Big\{ (t,x,y,z, a, b, \dot x,\dot z)\in\R^8\,\Big|\,1 + \dot x^2 + \dot z^2 + \5\,(2\,a - e^{-\sqrt{2}\,\omega\,x}\, b)^2 \, = \, a^2 \Big\} . $$ 
And the particular phase subspace $\,\EE_0\,$ distinguished in Remark \ref{rem.ptpart} can be written :   
$$ \EE_0 = \EE\cap \Big\{ a^2 = 1+ \dot z^2\,\,;\; 2\, a = e^{-\sqrt{2}\,\omega\,x}\, b\,\,;\; \dot x = 0\Big\} . $$

\brem \label{rem.pseudon} \  {\rm  
We see in particular that \ $a_s^2\ge 1\,$ and that \ $b_s^2 \ge 2\,e^{2\sqrt{2}\,\omega\,x_s}$, for any proper time $\,s\ge 0\,$. 
Therefore, \ $(a_s)$ and  $(b_s)$ almost surely never vanish. Moreover, they must have the same sign, since $ (00')$ implies \ ${\ds \Big| e^{-\sqrt{2}\,\omega\,x_s}\, {b_s\over a_s} -2\Big| \le \sqrt{2}\,}\,$ and then \  ${\ds e^{-\sqrt{2}\,\omega\,x_s}\, {b_s\over a_s} \ge 2-\sqrt{2}\,}$. This implies also \   ${\ds \Big| e^{\sqrt{2}\,\omega\,x_s}\, {a_s\over b_s} -1\Big| \le {1\over \sqrt{2}}}\,$. 
}\erem 

\subsection{Study of the one-dimensional sub-diffusions  $(a_s)$ and  $(\dot z_s)$} \label{SubDiffaz}\indf 
   These two one-dimensional sub-diffusions are easily handled. 
\blem \label{lem.a} \  There exist two standard real Brownian motions $(w_s)$ and $(w'_s)$, and two almost surely converging processes $(\eta_s)$ and $(\eta'_s)$, such that we have : 
$$ |a_s| = \exp\Big[ \sigma^2 s + \sigma\, w_s + \eta_s\Big] \quad \hbox{ for any proper time } \,s\ge 0\,, $$ 
and 
$$ |\dot z_s| = \exp\Big[ \sigma^2 s + \sigma\, w'_s + \eta'_s\Big]  \quad \hbox{ for any sufficiently large proper time }\, s \, . $$ 
\elem  
\ub{Proof} \quad  The stochastic differential equations satisfied by  $(a_s)$ and  $(\dot z_s)$ are respectively : 
$$ d a_s\, = \,  {\ts{3\,\sigma^2\over 2}}\, a_s\,ds + \sigma \, \sqrt{a_s^2-1}\, dw_s \, , \quad \hbox{
and }\quad d \dot z_s\, = \, {\ts{3\,\sigma^2\over 2}}\, \dot z_s\,ds + \sigma\,\sqrt{\dot z_s^2+1}\, dw'_s\, , $$ 
for two standard real Brownian motions $(w_s)$ and $(w'_s)$. \par 
   These equations are solved as follows ; \  we have real Brownian motions $(W_u)$ and $(W'_u)$ such that : 
$$ a_s = F \Big( W\Big[ \inf\Big\{ u\,\Big|\, \int_0^u (W_v^2-1)\2\, dv > \sigma\, s\Big\}\Big]\Big) , $$ 
and
$$ \dot z_s = G \Big( W'\Big[ \inf\Big\{ u\,\Big|\, \int_0^u (1-|W_v'|^2)\2\, dv > \sigma\, s\Big\}\Big]\Big) , $$ 
with \ $F(W) := {W\over  \sqrt{W^2-1}}\, $ and  \ $G(W') := {W'\over  \sqrt{1-|W'|^2}}\, $. \par 

   Clearly, as $\,u\,$ increases to the hitting time of 1 by $\,W$, then $\,{\ds\int_0^u (W_v^2-1)\2\, dv }\,$ increases to infinity, $\,W[u]$ goes to 1, and $\,F(W[u])$ goes to infinity, showing that $\,a_s\,$ goes almost surely to infinity with $\,s\,$. The same reasoning holds for $(\dot z_s)$ (except that $\,|W'_0|\,$ must be smaller than 1, while $\,|W_0|\,$ must be larger than 1), so that, in the same way, $\,|\dot z_s|\,$ goes almost surely to infinity with $\,s\,$.  \  The invariant measure of the diffusion $(a_s)$ is \ $1_{\{ |a| > 1\}} \sqrt{a^2-1}\, da\,$, and the invariant measure of the diffusion $(\dot z_s)$ is \ $\sqrt{z^2+1}\, dz\,$. \par  

Then  we have almost surely (since $\,|\dot z_s|\to\ii$), for any sufficiently large proper time $\,s\,$ : 
$$ d \log |a_s| = (1+ {\ts{1\over 2\,a_s^2}})\sigma^2 ds + \sigma \sqrt{1-a_s\2}\, dw_s \, ,  \hbox{
and }\; d \log |\dot z_s| = (1- {\ts{1\over 2\,\dot z_s^2}})\sigma^2 ds + \sigma\sqrt{1+\dot z_s\2}\, dw'_s\, . $$ 
Whence, for real Brownian motions $\,\tilde w\,,\,\tilde w'\,$ and for sufficiently large proper times $\,s_0, s$ : 
$$ \log |a_s| = \log |a_{0}| +\sigma^2 s + \sigma\, w_s + \int_0^s{\sigma^2 du\over 2\,a_u^2} - \tilde w\Bigg[\int_0^s{\sigma^2 \,a_u\2\,du\over 1+ \sqrt{1-a_u\2}}\Bigg] =  \sigma^2 s + o(s) > \sigma^2 s/2\, , $$ 
and similarly : 
$$ \log {|\dot z_s|\over |\dot z_{s_0}|} = \sigma^2 (s-s_0) + \sigma (w'_s-w'_{s_0}) - \int_{s_0}^s{{\ts\sigma^2} du\over 2\, \dot z_u^2} + \tilde w'\Bigg[\!\int_{s_0}^s\! {\sigma^2\,  \dot z_u\2\,du\over 1+\! \sqrt{1\! +\dot z_u\2}}\Bigg] \! =  \sigma^2 s + o(s)  > \sigma^2 s/2\, . $$ 
This implies the convergence of the integrals in the above formulas, and then the result. 
$\;\diamond$  

\bcor \label{cor.compz} \  For any sufficiently large proper time $\,s\,$ we have :  \ 
$ |z_s|\, = \,e^{\sigma^2 s + o(s^{5/9})}\, $.  
\ecor 

   We have also the following lower control, which we shall use later.  
\blem \label{lem.z} \  For any $\,A>\sqrt{3}\,$, we have \quad 
${\ds \P\Big[ (\exists\, s >0)\; |\dot z_s| \le A\,e^{\sigma^2 s/2}\Big/\,|\dot z_0| \ge A^2\Big]  < 1/\sqrt{A} \,}$, \  and \quad  ${\ds \P\Big[ (\exists\, s >0)\; |a_s| \le A\,e^{\sigma^2 s/2}\Big/\,|a_0| \ge A^2\Big]  < 1/A \,}$. 
\elem  
\ub{Proof} \quad  Fix $\,A>\sqrt{3}\,$ and $\,|\dot z_0| \ge A^2$. The stochastic differential equation satisfied by  $(\log|\dot z_s|)$, already written in the proof of Lemma \ref{lem.a}, is equivalent to : 
$$ d \log \Big|\dot z_s\,e^{-\sigma^2 s/2}\Big| =\, {\ts{\sigma^2\over 2}}\,(1- \dot z_s\2)\, ds + \sigma\sqrt{1+\dot z_s\2}\, dw'_s\, . $$ 
Let us apply the comparison theorem (see for example ([I-W], Theorem 4.1)) :  setting \par \smallskip  $T_A^z := \inf\{ s\,|\, |\dot z_s| = A\,e^{\sigma^2 s/2}\}$ \   and \  $ \log r_s := \log A^2 + {\ts{\sigma^2\over 2}}\,(1- A\2)\, s + \sigma\sqrt{1+A\2}\, w'_s\,$, \par\smallskip   \noindent  
we have : \qquad   
${\ds  \inf_{0\le s\le T_A^z}\, \Big|\dot z_s\,e^{-\sigma^2 s/2}\Big| \ge \inf_{0\le s\le T_A^z}\, r_s\,}$, \ whence 
$$ \P[ T_A^z <\ii] \le \P[ \;\log r_s \; \hbox{\it hits } \log A ] = \P[ \; w'_s-\,{\ts{1-A\2\over 1+A\2}}\, s/2 \;\, \hbox{\it hits } \log A ] = A^{-{1-A\2\over 1+A\2}} < 1/\sqrt{A} \, . $$ 
Similarly, for $\,|a_0| \ge A^2\,$ and \  $T_A^a := \inf\{ u\,|\, |a_u| = A\,e^{\sigma^2 s/2}\}$, \   since 
$$ d  \Big|a_s\,e^{-\sigma^2 s/2}\Big| = \, {\ts{\sigma^2\over 2}}\,a_s\2\, ds + \sigma \sqrt{1-a_s\2}\, dw_s \, , $$ 
we get : 
$$ \P[ T_A^a <\ii] \le \P[ \; w_s-s/2 \;\, \hbox{\it hits } \log A ] = 1/A \, . \;\;\diamond $$

\subsection{Study of the two-dimensional sub-diffusion $(a_s,\dot z_s)$} \label{SubDiff2} \indf 
  Recall from Section \ref{sec.RedDim} that we have : \  $\,a_s^2- 1\ge \dot z_s^2\,$,  
$$ d a_s\, =\,  {\ts{3\,\sigma^2\over 2}}\, a_s\,ds + \sigma \, dM^a_s  \quad  \hbox{
and }\quad   d \dot z_s\, =\, {\ts{3\,\sigma^2\over 2}}\, \dot z_s\,ds + \sigma \, dM^z_s \, , $$ 
where the $\,\R^2$-valued martingale \   $\check M_s := (M^a_s, M^z_s)$ has quadratic covariation matrix : 
$$ ((\check K^{ij}_s)) =  \pmatrix{ a_s^2 -1 &  a_s\,\dot z_s \cr 
a_s\,\dot z_s & \dot z_s^2+1 \cr } . $$ 

\if{ 
\bpro  \label{pro.az} \  The process $\,\dot z_s/a_s\,$ converges almost surely, toward some random limit $\,\ell\,$ such that $\,0<|\ell | \le 1\,$. 
\epro 
\ub{Proof} \quad    Using It\^o's Formula and the above expressions for $\,d a_s \,$ and  $\,d \dot z_s \,$, we get : 
$$ d\Bigg[ {\dot z_s\over \sqrt{a_s^2-1}}\Bigg] =\,  {d\dot z_s\over \sqrt{a_s^2-1}} - {\dot z_s\,a_s\,d a_s\over (a_s^2-1)^{3/2}} - {a_s\,\langle d a_s, d\dot z_s\rangle\over (a_s^2-1)^{3/2}} + {(a_s^2+\5) \,\dot z_s \,\langle d a_s, d a_s\rangle\over (a_s^2-1)^{5/2}} $$ 
$$ =\,  {-\,\sigma^2\, \dot z_s\over (a_s^2-1)^{3/2}}\, ds + {\sigma\over  \sqrt{a_s^2-1}}\, \sqrt{ 1- {\dot z_s^2\over a_s^2-1}}\, d\check w_s  \, , $$ 
and then (using Lemma \ref{lem.a}), for large enough proper time $\,s$ :  
$$ d\log \Bigg[ {|\dot z_s|\over \sqrt{a_s^2-1}}\Bigg] =\, {-\,\sigma^2\over a_s^2-1}\, ds + {\sigma\over \dot z_s}\, \sqrt{ 1- {\dot z_s^2\over a_s^2-1}}\, d\check w_s -  { \sigma^2\over 2\, \dot z_s^2}\times \bigg[ 1- {\dot z_s^2\over a_s^2-1}\bigg] ds \, . $$ 
Therefore, for almost any trajectory, fixing some large enough proper time $\,s_0\,$, for any  $\,s\ge s_0\,$ we have (for some real Brownian motion $\,\check W$) : 
$$ \log \Bigg[ {|\dot z_s|\over \sqrt{a_s^2-1}}\Bigg] =\, \log \Bigg[ {|\dot z_{s_0}|\over \sqrt{a_{s_0}^2-1}}\Bigg] - {\sigma^2\over 2} \bigg[ \int_{s_0}^s {du\over a_u^2-1} + \int_{s_0}^s {du\over \dot z_u^2}\bigg] + \sigma\, \check W\bigg[ \int_{s_0}^s {du\over \dot z_u^2} -\int_{s_0}^s {du\over a_u^2-1} \bigg] , $$ 
which almost surely converges in $\,\R\,$, by Lemma \ref{lem.a}. By Lemma \ref{lem.a} again, this implies clearly the almost sure convergence of $\,\dot z_s/a_s\,$ toward some random limit $\,\ell\in\R^*$. Finally the unit pseudo-norm relation $(00')$ implies $\,a_s^2- 1\ge \dot z_s^2$, and then $\,\ell^2\le 1\,$. $\;\diamond$  
}\fi 
   We get easily an asymptotic variable of the sub-diffusion $(a_s,\dot z_s)$. 
\bpro  \label{pro.az} \  The process $\,(\dot z_s/a_s)\,$ converges almost surely, toward some random limit $\,\ell\,$ such that $\,0<|\ell | \le 1\,$.  
\epro 
\ub{Proof} \quad  Using It\^o's Formula and the above expressions for $\,d a_s \,$ and  $\,d \dot z_s \,$, we get : 
$$ d \Big[{\dot z_s\over a_s}\Big]\, =\, {d \dot z_s\over a_s} - {\dot z_s\, d a_s\over a_s^2} + {\dot z_s\, \langle d a_s\rangle \over a_s^3} - {\langle d a_s, d \dot z_s\rangle \over a_s^2}  =\,  {\sigma\over a_s}\Big[ dM^z_s- {\dot z_s\over a_s}\,dM^a_s\Big] -  {\sigma^2\, \dot z_s\over a_s^{3}}\, ds \, , $$ 
with 
$$ \Big\langle a_s\1\Big[ dM^z_s- {\dot z_s\over a_s}\,dM^a_s\Big]\Big\rangle\, =\, {1-|\dot z_s/a_s|^2\over a_s^2}\, ds \, . $$  
Hence, we have some real Brownian motion $\,\check W\,$ such that almost surely, for any $\,s\ge 0$ :  
$$ {\dot z_s\over a_s}\, =\, {\dot z_0\over a_0} - \sigma^2 \int_{0}^s  {\dot z_u\over a_u^{3}}\, du + \sigma \, \check W\bigg[ \int_{0}^s {1-|\dot z_u/a_u|^2\over a_u^2}\, du \bigg] , $$ 
which almost surely converges toward some random limit $\,\ell\in\R$, by Lemma \ref{lem.a}. \  The unit pseudo-norm relation $(00')$ implies $\,a_s^2- 1\ge \dot z_s^2$, and then $\,\ell^2\le 1 $. \par 
   Similarly, using Lemma \ref{lem.a} again, for any large enough proper time $\,s\,$ we have :  
$$ d \Big[{a_s\over \dot z_s}\Big]\, =\, {d a_s\over \dot z_s} - {a_s\, d \dot z_s\over \dot z_s^2} + {a_s\, \langle d \dot z_s\rangle \over \dot z_s^3} - {\langle d a_s, d \dot z_s\rangle \over \dot z_s^2}  =\,  {\sigma\over \dot z_s}\Big[ dM^a_s- {a_s\over \dot z_s}\,dM^z_s\Big] + {\sigma^2\, a_s\over \dot z_s^{3}}\, ds \, , $$ 
with 
$$ \Big\langle \dot z_s\1\Big[ dM^a_s- {a_s\over \dot z_s}\,dM^z_s\Big]\Big\rangle\, =\, {|a_s/\dot z_s|^2 -1 \over \dot z_s^2}\, ds \, , $$  
whence the almost sure convergence of $(a_s/\dot z_s)$, which proves that $\,\ell\not= 0\,$ almost surely. $\;\diamond$  
\par \medskip 

   The following statement ensures that the range of possible limits $\,\ell\,$ in Proposition \ref{pro.az}, is the whole $[-1,0[\,\cup\, ]0,1]$. This provides a continuum of non-trivial bounded harmonic functions for the relativistic operator $\LL\,$. 
\bpro  \label{pro.a/z} \  For any real $\,\ell_0\,$ such that $\,0< |\ell_0| \le 1\,$, and for any $\,\e>0\,$, we have \quad  ${\ds\,\P \Big[\ell_0-\e < \ell = \lim_{s\to\ii}\,{\dot z_s\over a_s}\, < \ell_0+\e \Big] > 1-\e\,}$, provided \  ${\dot z_0/a_0}\,$ is close enough from $\,\ell_0\,$ and $\,|a_0|\,$ is large enough.
\epro 
\ub{Proof} \quad  Fix $\,A> 9\,$, $\,|a_0| > A^2$ and  $\,|\dot z_0| \ge A^2$, such that \  ${\dot z_0/a_0}\,$ is close to  $\,\ell_0\,$ (precisely, we demand $\,|\log({\dot z_0\over a_0\,\ell_0})| <A\2$), and consider the event : $$ \AA\, := \, \Big\{\,a_s^2 > 1+ A^2 e^{\sigma^2 s} \;\hbox{ and }\;\; \dot z_s^2 > A^2 e^{\sigma^2 s} \quad \hbox{ for all } s \ge 0 \,\Big\} . $$ 
By Lemma \ref{lem.z} we have : \  $\P(\AA) > \, 1-2/\sqrt{A}\,$. \  Now, on $ \AA\,$ we have : 
$$ \int_{0}^\ii {du\over a_u^2-1} + \int_{0}^\ii {du\over \dot z_u^2} \,\le \, 2\sigma\2\,A\2 \quad \hbox{and} \quad   \int_{0}^\ii {du\over \dot z_u^2} \,\le \, \sigma\2\,A\2\, . $$ 
Hence, we see from the expression giving \ ${\ds  \log \Bigg[ {|\dot z_s|\over \sqrt{a_s^2-1}}\Bigg]}$, displayed in the proof of Proposition \ref{pro.az}, that we have on $ \AA$ :  
$$ |\log(\ell/\ell_0)| \le 2\, A\2 + \sigma\,\max\{ |\check W_s|\,|\, 0\le s\le \sigma\2 A\2\}\, . $$ 
Finally, as \  
$$ \P[\sigma\,\max\{ |\check W_s|\,|\, 0\le s\le \sigma\2 A\2\}> A^{-1/2}] \le 2\, \P[\max\{ \check W_s\,|\, 0\le s\le A\2\}> A^{-1/2}]  $$ 
$$ = 2\, \P[\,|\check W_{A\2}| > A^{-1/2}] = 4\, \P[\check W_{1} > \sqrt{A}\,] < e^{-A/2}\, ,  $$ 
we obtain : \ ${\ds \P[ \, |\log(\ell/\ell_0)| \le 2\, A\2 +A^{-1/2}] > 1-2/\sqrt{A}\, - e^{-A/2}\,}$. $\;\diamond$ 

\bpro  \label{pro.a-z} \  The process $\,\lambda_s := \argsh \Big(\sqrt{a_s^2-\dot z_s^2-1}\,\Big)\,$ solves the following stochastic differential equation : \quad   
${\ds  d \lambda_s = \sigma\,  d\tilde w_s + \sigma^2\,\coth(2\lambda_s)\, ds \, , }$\  
for some standard real Brownian motion $\,\tilde w\,$, and then is a real diffusion which goes almost surely to infinity as $\,s\to\ii\,$. Moreover, the process \  $\;{ \tilde\eta_s:= \lambda_s - \sigma \tilde w_s- \s^2 s}\,$ converges almost surely  (in $\R$) as $\,s\to\ii\,$. 
\epro 
\ub{Proof} \quad    Using It\^o's Formula and the expressions for $\,d a_s \,$ and  $\,d \dot z_s \,$, we get : \par
\vbox{ 
$$ d\lambda_s\, =\,  {2\, a_s\over \sh(2\lambda_s)}\,\Big[ {\ts{3\,\sigma^2\over 2}}\, a_s\,ds + \sigma \, dM^a_s\Big] - {2\, \dot z_s\over \sh(2\lambda_s)}\,\Big[ {\ts{3\,\sigma^2\over 2}}\, \dot z_s\,ds + \sigma \, dM^z_s\Big] + {a_s^2\dot z_s^2\, \ch(2\lambda_s)\over \ch^3\lambda_s\, \sh^3\lambda_s}\, \sigma^2  ds  $$ 
$$ - { \sigma^2\over  \sh(2\lambda_s)}\,\Big[  {a_s^2\over \ch^2\lambda_s} + {\dot z_s^2+1\over \sh^2\lambda_s}\Big] (a_s^2-1)\, ds - { \sigma^2\over  \sh(2\lambda_s)}\,\Big[  {\dot z_s^2\over \ch^2\lambda_s} + {a_s^2-1\over \sh^2\lambda_s}\Big] (\dot z_s^2+1)\, ds  $$ } \par
\smallskip  \vbox{ 
$$ =\, {2\,\sigma \over \sh(2\lambda_s)} \Big[ a_s\,dM^a_s - \dot z_s\,dM^z_s\Big] + {3\,\sigma^2\, \ch^2\lambda_s \over \sh(2\lambda_s)}\, ds + {a_s^2\dot z_s^2\, \ch(2\lambda_s)\over \ch^3\lambda_s\, \sh^3\lambda_s}\, \sigma^2  ds  $$
$$ -\, { \sigma^2\over  \sh(2\lambda_s)}\,\Big[  {a_s^4+\dot z_s^4 - \ch^2\lambda_s\over \ch^2\lambda_s} + 2\,{a_s^2\,\dot z_s^2 + \sh^2\lambda_s\over \sh^2\lambda_s}  \Big] \, ds $$ } 
\par \vspace{-2mm}  \vbox{ 
$$ = \, {2\,\sigma \over \sh(2\lambda_s)}\Big\langle (a_s^2-\dot z_s^2)(a_s^2-\dot z_s^2-1)\Big\rangle^{1/2} d\tilde w_s +  { \sigma^2\over  \sh(2\lambda_s)}\,\Big[ 3\,\ch^2\lambda_s + 2\,
{a_s^2\dot z_s^2\over \ch^2\lambda_s} -{a_s^4+\dot z_s^4\over \ch^2\lambda_s} -1\Big] \, ds $$ } 
\vspace{-6mm} 
$$ =\, \sigma\,  d\tilde w_s + \sigma^2\,\coth(2\lambda_s)\, ds \, . $$

   Since $\,\coth(2\lambda_s) > 1\,$, the comparison theorem ensures that we have almost surely : 
$$ \lambda_s \ge \lambda_0 + \sigma\,\tilde w_s + \sigma^2\,s \lra +\ii\, . $$ 
Moreover, we have almost surely for large enough $\,s_0\,$ and for $\,s\ge s_0$ : \quad  
${\ds \lambda_s \ge \sigma^2 s/2\,}$. \  Hence, we deduce that \quad ${\ds \tilde\eta_s = \tilde\eta_{s_0}+ 2\sigma^2\int_{s_0}^s{du\over e^{4 \lambda_u}-1}\,}\;$ converges almost surely.   $\;\diamond $ 

\brem \label{rem.lambda} \   {\rm The equation satisfied by $(\lambda_s)$ can be precisely solved as follows, provided $\,\lambda_0>0\,$, using some real Brownian motion $\,\beta $ started from $\,\beta_0=  \5\, \log [ \coth\lambda_0]$ : 
$$ \lambda_s = \5\, \log\bigg[ \coth\Big( \beta \Big[ \inf\Big\{ u\,\Big|\,\int_0^u\sh\2(2\beta_v)\, dv = \sigma^2 s\Big\}\Big]\Big)\bigg] . $$ 
This implies that we have almost surely : \  $\lambda_s>0\,$ for any $\,s >0 $ : \  the state subspace $\,\EE_0\,$ of Remark \ref{rem.ptpart} is polar for the relativistic diffusion. Recall from Section \ref{sec.RD} (before Section \ref{sec.RedDim}) that it is also instantly unstable. Hence we can always restrict the state space $\,\EE\,$ of the relativistic diffusion to $\,\EE\,{\ss\sm}\,\EE_0\,$. 
}\erem 

      By the unit pseudo-norm relation $(00')$ and by Remark \ref{rem.lambda}, we have  almost surely for any $\,s>0$ :  $\,a_s^2- 1> \dot z_s^2$. Therefore there exist two independent standard real Brownian motions $(w_s)$ and $(\check w_s)$  such that : 
$$ d a_s = {\ts{3\,\sigma^2\over 2}}\, a_s\,ds + \sigma\sqrt{a_s^2-1}\, dw_s \;\,  \hbox{
and }\;\, d \dot z_s = {\ts{3\,\sigma^2\over 2}}\, \dot z_s\,ds + \sigma\,{a_s\, \dot z_s\over  \sqrt{a_s^2-1}}\, dw_s + \sigma\sqrt{ 1- {\dot z_s^2\over a_s^2-1}}\, d\check w_s  \, . $$ 

   We shall need to know that in fact $\,|\ell | < 1\,$ almost surely. \   Let us set : 
$$ (11) \qquad  A_s := \sqrt{1 - {1+\dot z_s^2\over a_s^2}}\,\, = \, \sqrt{\Big[{\dot x_s\over a_s}\Big]^2 + \5\, \Big[ 2 - e^{-\sqrt{2}\,\omega\,x_s}\, {b_s\over a_s} \Big]^2}\, 
= \bigg|{\sh \lambda_s\over a_s}\bigg| \, . $$ 
Note that the phase subspace $\,\EE_0\,$ (of Remarks \ref{rem.ptpart} and \ref{rem.lambda}) is precisely :   \  $\,{\ds \EE_0 = \EE\cap \{ A = 0 \} . }$
   
\bpro \label{pro.ell} \   The random limit $\,\ell =\lim_{s\to\ii}(\dot z_s/a_s)\,$ of Proposition \ref{pro.az} satisfies almost surely : $\,0<|\ell | < 1\,$. 
\epro 
\ub{Proof} \   By Proposition \ref{pro.az} and Notation (11), we have almost surely : 
\  $ \log(1-\ell^2) = 2\, \log A_\ii \,$. \  On the other hand, using It\^o's Formula and Proposition \ref{pro.a-z}, we get : 
$$ d (\log A_s) = d (\log [\sh\lambda_s]) - d (\log |a_s|) = \coth\lambda_s \, d\lambda_s - { d\langle \lambda_s\rangle\over 2\, \sh^2\lambda_s} - \sigma^2(1+\5\,a_s\2) ds - {\sigma\over a_s}\, dM^a_s $$  
$$ =\, {\sigma} \bigg[\bigg( {a_s\over \sh^2\lambda_s} - {1\over a_s}\bigg) dM^a_s - {\dot z_s\over \sh^2\lambda_s}\,dM^z_s\bigg] + {\sigma^2\over 2} \bigg[ 2\,\coth\lambda_s\,\coth(2\lambda_s) - {1\over \sh^2\lambda_s} - 2 - {1\over a_s^2} \bigg] ds $$ 
$$ =\, {\sigma\over a_s\,\sh^2\lambda_s} \Big[ (\dot z^2_s+1)\, dM^a_s - a_s\,\dot z_s\,dM^z_s\Big] - {\sigma^2\over 2\, a_s^2} \, ds \, , $$ 
whence, for some real Brownian motion $\,\check B$ :  
$$  \log(1-\ell^2) = \, 2\, \log A_0 + 2\,\sigma\, \check B\bigg[ \int_0^\ii {\dot z_s^2 +1\over a_s^2\,(a_s^2- \dot z_s^2 -1)}\, ds \bigg] - \sigma^2\int_0^\ii {ds\over a_s^2} \,, $$ 
which converges (in $\,\R$) almost surely, by Lemma \ref{lem.a}, Proposition \ref{pro.az} and Proposition \ref{pro.a-z}, showing that indeed $\,\ell^2< 1\,$ almost surely. $\;\diamond$  
\par\medskip

   We have furthermore the following. 
\bpro \label{lem.azn} \  \if{ $(i)$ \  We have almost surely, for any sufficiently large proper time $\,s$ : 
$$ |\dot z_s- \ell\,a_s|\, = \, \exp\Big[ \sigma^2 s + o(s^{5/9})\Big] \, ,\;\hbox{ and} \quad \Big| z_s- \ell\int_0^sa_u\, du \Big|\, = \, \exp\Big[ \sigma^2 s + o(s^{5/9})\Big]  . $$ \par 
$(ii)$ \}\fi   The law of the random limit $\,\ell =\lim_{s\to\ii}\limits (\dot z_s/a_s)\,$ has no atom. 
\epro  
\ub{Proof} \quad  Fix any $\,\ell_0\in\,]-1,1[\,$, and set \   $\delta_s := \dot z_s- \ell_0\,a_s\,$. The stochastic differential equation satisfied by $(\delta_s)$ is easily seen to be : 
$$ d \delta_s\, = \,  {\ts{3\,\sigma^2\over 2}}\, \delta_s\,ds + \sigma \, \sqrt{\delta_s^2+1-\ell_0^2}\, \, d\beta_s \, , $$ 
for some standard real Brownian motion $(\beta_s)$. This diffusion equation can be solved as follows : \  we have a real Brownian motion $(W_u)$ (started from $\,W_0\in \,]{-1\over 1-\ell_0^2},{1\over 1-\ell_0^2}[\,$) such that : 
$$ \delta_s = F \Big( W\Big[ \inf\Big\{ u\,\Big|\, \int_0^u {(1-\ell_0^2)\, dv \over 1- (1-\ell_0^2)^{2}\,W^2_v}\, > \sigma\, s\Big\}\Big]\Big) , $$ 
with \ $F(W) := {(1-\ell_0^2)^{3/2}\,W\over  \sqrt{1- (1-\ell_0^2)^{2}\,W^2}}\, $. \quad   As $\,u\,$ increases to the hitting time of $\,\pm(1-\ell_0^2)\1\,$ by $\,W$, then $\,{\ds\int_0^u {(1-\ell_0^2)\, dv \over 1- (1-\ell_0^2)^{2}\,W^2_v}\, }\,$ increases to infinity, $\,W_u\,$ goes to $\,\pm(1-\ell_0^2)\1$, and $\,F(W_u)$ goes to $\,\pm\ii\,$, showing that $\,|\delta_s|\,$ goes almost surely to infinity with $\,s\,$.  \   (The invariant measure of the diffusion $(\delta_s)$ is  \ $\sqrt{\delta^2+1-\ell_0^2}\,\, d\delta\,$.) \par  

Then  we have almost surely, for any sufficiently large proper time $\,s\,$ : 
$$ d \log |\delta_s| = (1- {\ts{1-\ell_0^2\over 2\,\delta_s^2}})\,\sigma^2 ds \pm \sigma\sqrt{1+{\ts{1-\ell_0^2\over \delta_s^2}}}\, d\beta_s\, . $$ 
Whence, for real Brownian motions $\,w,\tilde w\,$ and for sufficiently large proper times $\,s_0, s$ : 
$$ \log {|\delta_s|\over |\delta_{s_0}|} = \sigma^2 (s-s_0) + \sigma\, (w_s-w_{s_0}) - 
{\ts {\sigma^2(1-\ell_0^2)\over 2}}\! \int_{s_0}^s{du\over \delta_u^2} + {\ss \sigma^2(1-\ell_0^2)}\,  \tilde w\bigg[\!\int_{s_0}^s\!{du\over \delta_u^2}\bigg]\! =  \sigma^2 s + o(s)  > \sigma^2 s/2\, . $$ 
This implies the convergence of the integrals in the above formula, and then the existence of a standard real Brownian motion $(w_s^{\ell_0})$ and of an almost surely converging process $(\eta^{\ell_0}_s)$, such that almost surely, for any sufficiently large proper time $\,s\,$ we have : 
$$ |\delta_s| = |\dot z_s- \ell_0\,a_s| = \exp\Big[ \sigma^2 s + \sigma\, w^{\ell_0}_s + \eta^{\ell_0}_s\Big] = \exp\Big[ \sigma^2 s + o(s^{5/9})\Big] \, . $$ 

   For the same $\,\ell_0\,$ and $(\delta_s)$ as above, we get as in the proof of Proposition \ref{pro.az} : 
$$ \log {|\delta_s/a_s|\over |\delta_{s_0}/a_{s_0}|} = - \sigma^2 \int_{s_0}^s {\dot z_u\,du\over a_u^2\,\delta_u} - {\sigma^2\over 2} \int_{s_0}^s {1-(\dot z_u/a_u)^2\,du\over \delta_u^2}
+ {\sigma}\,  \check W\bigg[\!\int_{s_0}^s{1-(\dot z_u/a_u)^2\,du\over \delta_u^2}\bigg] , $$ 
almost surely for any sufficiently large $\,s_0, s\,$. Using the above, this shows the almost sure convergence of \ $\log\Big|{\dot z_s\over a_s}-\ell_0\Big|$, hence by Proposition \ref{pro.az}, that indeed $\,\P[\ell =\ell_0] =0\,$.   $\;\diamond$  
\par\medskip 

   The following statement implies that the asymptotic $\,\sigma$-algebra of $(a_s,\dot z_s)$ is generated by the only variable $\,\ell\,$. 
\bpro \label{pro.az''} \  The process $(\dot z_s - \ell\,a_s)$ converges in law, toward some 
smooth law, but not in probability. Furthermore, we have \ $(\dot z_s - \ell\,a_s) = o(|a_s|^{2/9})$ almost surely. 
\epro 
\ub{Proof} \quad  Since \quad ${\ds \Big\langle{dM^a_s\over a_s}\Big\rangle = (1-a_s\2) ds \,}$, 
$$  \Big\langle{ dM^z_s- {\dot z_s\over a_s}\,dM^a_s}\Big\rangle = \Big(1-\Big|{\dot z_s\over a_s}\Big|^2\Big) ds\; , \; \hbox{ and }\quad \Big\langle{ {dM^a_s\over a_s}\, ,\,dM^z_s- {\dot z_s\over a_s}\,dM^a_s}\Big\rangle = {\dot z_{s}\over a_{s}^2}\,  ds \, , $$
there exist two independent standard real Brownian motions $\,w,w'\,$ such that : 
$$ {dM^a_s\over a_s} = \sqrt{1-a_s\2}\, dw_s \; , \; \hbox{ and }\quad dM^z_s- {\dot z_s\over a_s}\,dM^a_s = {\dot z_s\over a_s^2\sqrt{1-a_s\2}}\,dw_s + \sqrt{a_s^2-1-\dot z_s^2\over a_s^2-1}\, dw'_s\, . $$ 
Then, since
$$ d\log a_s = \sigma^2 (1+ {\ts{1\over 2\,a_s^2}}) ds + {\sigma\over a_s}\, dM^a_s\, ,  $$ 
we have for any $\,0\le s\le u$ (recall Lemma \ref{lem.a}) : 
$$ {a_s\over a_{s+u}} \,=\, \exp\bigg[ -\sigma^2 u - {\ts{\sigma^2\over 2}}\!\int_s^{s+u} {dv\over a_v^2} - \sigma\! \int_s^{s+u} \sqrt{1-a_v\2}\, dw_v \bigg] = e^{-\sigma^2 u - \sigma (w_{s+u} -w_s) + o(1)} . $$    
Hence, the expression for $\,d[\dot z_s/a_s]$ used for Proposition \ref{pro.az}, implies that almost surely  :   
$$ {\dot z_s} -\ell\,a_s \, =\, \sigma^2 a_s\! \int_s^\ii  {\dot z_u\over a_u^{3}}\, du - \sigma\,a_s \int_s^\ii a_u\1 \Big[ dM^z_u- {\dot z_u\over a_u}\,dM^a_u\Big]  $$ 
$$ =\, \sigma^2 a_s\! \int_s^\ii  \O(a_u^{-2})\, du  - \sigma\,a_s\! \int_s^\ii\! {\dot z_{u}\, dw_{u}\over a_{u}^3 \sqrt{1-a_{u} \2}}  - \sigma\,a_s\! \int_s^\ii\! \sqrt{1-{\dot z_{u}^2\over a_{u}^2-1}}\,  \,{dw'_{u}\over a_{u}} \, .  $$ 
Now, there exist Brownian Motions $\tilde W, \tilde W'$ (independent of $\,s$) such that : 
$$ \int_s^\ii \! {\dot z_{u}\, dw_{u}  \over a_{u}^3 \sqrt{1-a_{u} \2}} = \tilde W\bigg[\! \int_s^\ii \! {\dot z_{u}^2\, d{u}  \over a_{u}^6 (1-a_{u} \2)}\bigg] = o\bigg[\!\int_s^\ii \! {\dot z_{u}^2\, d{u}  \over a_{u}^6 (1-a_{u} \2)}\bigg]^{4/9} = o\bigg[\! \int_s^\ii \! a_{u}^{-4}\,{du}\bigg]^{4/9} $$ 
$$ =\, o\Big(|a_s|^{-16/9}\Big) \times \bigg[\! \int_0^\ii \Big|{a_s\over a_{s+u}}\Big|^{4} {du}\bigg]^{4/9} =\, o\Big(|a_s|^{-16/9}\Big) \times \bigg[\! \int_0^\ii e^{- 4 \sigma^2 u - 4\sigma (w_{s+u}-w_s)}  du\bigg]^{4/9} $$ 
$$ =\,  o\Big(|a_s|^{-4/3}\Big) \times \bigg[\! \int_0^\ii e^{- 4 \sigma^2 u - 4\sigma (w_{s+u}-w_s)- \sigma^2 s - \sigma w_s}  du\bigg]^{4/9} $$ 
$$ =\,  o\Big(|a_s|^{-4/3}\Big) \times \bigg(\!\int_0^\ii\! \exp\!\Big(- \Big[ 7\,u + (s+u) \Big( 1+ {\ts{8\,w_{s+u}\over \sigma\,(s+u)}} \Big) + s \Big(1-{\ts{6\,w_{s}\over \sigma\,s}} \Big)\Big] \sigma^2/ 2\Big) du \bigg)^{4/9} =  o\Big(|a_s|^{-4/3}\Big) , $$ 
and similarly : 
$$ \int_s^\ii\! \sqrt{1-{\dot z_{u}^2\over a_{u}^2-1}}\,\,{dw'_{u}\over a_{u}} = \tilde W'\bigg[\! \int_s^\ii \! \Big[{1-{\dot z_{u}^2\over a_{u}^2-1}}\Big] {d{u}  \over a_{u}^2}\bigg] =  o\bigg[\! \int_s^\ii \! a_{u}^{-2}\,{du}\bigg]^{4/9} $$ 
$$ =\, o\Big(|a_s|^{-8/9}\Big) \times \bigg[\! \int_0^\ii \Big|{a_s\over a_{s+u}}\Big|^{2} {du}\bigg]^{4/9} =\, o\Big(|a_s|^{-8/9}\Big) \times \bigg[\! \int_0^\ii e^{- 2 \sigma^2 u - 2\sigma (w_{s+u}-w_s)}  du\bigg]^{4/9} $$ 
$$ =\,  o\Big(|a_s|^{-2/3}\Big) \times \bigg(\!\int_0^\ii\! \exp\!\Big(- \Big[ 7\,u + (s+u) \Big( 1+ {\ts{8\,w_{s+u}\over \sigma\,(s+u)}} \Big) + s \Big(1-{\ts{6\,w_{s}\over \sigma\,s}} \Big)\Big] \sigma^2/ 4\Big) du \bigg)^{4/9} =  o\Big(|a_s|^{-7/9}\Big) . $$ 
This shows also that \  ${\ds \int_s^\ii  \O(a_u^{-2})\, du = o\Big(|a_s|^{-7/4}\Big)}$. \par \smallskip 

   So far, we have shown that almost surely, as $\,s\to\ii$ : 
$$ {\dot z_s} -\ell\,a_s\,  = \, o\Big(|a_s|^{-1/3}\Big) 
- \sigma\,a_s\! \int_s^\ii\! \sqrt{1-{\dot z_{u}^2\over a_{u}^2-1}}\,  \,{dw'_{u}\over a_{u}}\, =\,  o\Big(|a_s|^{2/9}\Big) . $$ 

    Considering then a standard real Brownian motion $\,w''$, independent from $(a_s,\dot z_s)$, we have : 
$$ \Big\langle \int_s^\ii\! {a_s\over a_u} \bigg[{\ts\sqrt{1-{\dot z_{u}^2\over a_{u}^2-1}}}\,-\sqrt{1-\ell^2}\,\bigg] dw''_{u}\, \Big\rangle $$
$$  = \int_0^\ii\! e^{-2\sigma^2 u - 2\sigma (w_{s+u} -w_s)+ o_s(1)} \bigg[{\ts\sqrt{1-{\dot z_{s+u}^2\over a_{s+u}^2-1}}}\,-\sqrt{1-\ell^2}\,\bigg]^2 du $$ 
$$ \le \sqrt{\int_0^\ii\! e^{-2\sigma^2 u} \bigg[{\ts\sqrt{1-{\dot z_{s+u}^2\over a_{s+u}^2-1}}}\,-\sqrt{1-\ell^2}\,\bigg]^4 du} \times 
\sqrt{\int_0^\ii\! e^{-2\sigma^2 u - 4\sigma (w_{s+u} -w_s)+ o_s(1)}\, du} \,\, , $$ 
which goes to 0 in probability as $\,s\to\ii\,$.  Hence, $({\dot z_s} -\ell\,a_s)$ behaves in probability as : 
$$ -\sigma \sqrt{1-\ell^2} \int_s^\ii\! e^{-\sigma^2 (u-s) - \sigma (w_{u} -w_s)} \, dw''_{u}\, ,  $$ 
and then converges in law, toward the law of : 
 $$ \sqrt{1-\ell^2} \int_0^\ii\! e^{- u - w_{u}} \, dw''_{u}\, \equiv\, \sqrt{(1-\ell^2)\! \int_0^\ii\! e^{-2\, u - 2\, w_{u}} \, du} \times N\,  ,  $$ 
$N\,$ denoting a $\,\NN(0,1)$ Gaussian variable, independent from $(\ell, w)$. \par \smallskip   
 
        On the other hand, the above proves also that, as $\,t\ge s\to\ii\,$,  \  
$ ({\dot z_t} -\ell\,a_t) - ({\dot z_s} -\ell\,a_s) $ \   behaves in probability as : 
$$ \sigma \sqrt{1-\ell^2}\, \bigg[\int_s^\ii\! e^{-\sigma^2 (u-s) - \sigma (w_{u} -w_s)} \, dw''_{u} - \int_t^\ii\! e^{-\sigma^2 (u-t) - \sigma (w_{u} -w_t)} \, dw''_{u}\bigg]\, ,  $$ 
which diverges in probability :  for small enough $\,\e>0\,$ and large enough $(t-s)$, we have 
${\ds \,\P\bigg[\Big|\int_t^\ii\! e^{-\sigma^2 (u-s) - \sigma (w_{u} -w_s)} \, dw''_{u}\Big| <\e \bigg] >1-\e\,}$, and by independence : 
$$ \P\bigg[\Big|\int_s^t\! e^{-\sigma^2 (u-s) - \sigma (w_{u} -w_s)} \, dw''_{u} - \int_t^\ii\! e^{-\sigma^2 (u-t) - \sigma (w_{u} -w_t)} \, dw''_{u}\Big| > 2\e \bigg] >1-\e \, . \;\;\diamond $$ 

\if{
   $(ii)$ \  By Lemma \ref{lem.a} and Corollary \ref{cor.compz}, the above expression for $\,\dot z_s/a_s\,$ implies : 
$$ {\dot z_s\over a_s} -\ell \, =\, \sigma^2 \int_s^\ii  {\dot z_u\over a_u^{3}}\,  du + \sigma \, \check W\bigg[ \int_s^\ii {1-|\dot z_u/a_u|^2\over a_u^2}\, du \bigg]  $$ 
$$ =\, \int_s^\ii  e^{-2\sigma^2 u + o(u^{5/9})}\, du + \sigma \, \check W\bigg[ \int_s^\ii e^{-2\sigma^2 u + o(u^{5/9})}\, du \bigg] = e^{-\sigma^2 s + o(s^{5/9})}\, , $$ 
hence \  ${\ds \int^\ii \Big|{\dot z_s\over a_s} -\ell\Big|\, ds \, <\ii \,}$ almost surely. \par \smallskip 
   Then we have successively :
$$ d \bigg[ {\log |a_s|\over a_s}\bigg] =\, {1-\log |a_s|\over a_s^2}\,d a_s + {2\log |a_s| -3\over 2\,a_s^3}\,\langle d a_s \rangle $$
$$ =\, {\sigma\over a_s^2}\,(1-\log |a_s|)\, dM^a_s - {\sigma^2 \log |a_s|\over 2\,a_s}\, ds - {\sigma^2\over 2\,a_s^3}\,(2\log |a_s| -3)\, ds\, ; $$ 
\vbox{
$$ d \bigg[ \dot z_s\times {\log |a_s|\over a_s}\bigg] =\, {\sigma\, \dot z_s\over a_s^2}\,(1-\log |a_s|)\, dM^a_s - {\sigma^2\,\dot z_s \log |a_s|\over 2\,a_s}\, ds - {\sigma^2\,\dot z_s\over 2\,a_s^3}\,(2\log |a_s| -3)\, ds \qquad $$ 
$$ \qquad +\, {3\,\sigma^2\,\dot z_s\,\log |a_s|\over 2\,a_s}\, ds + \sigma\,{\log |a_s|\over a_s}\, dM^z_s +  {\sigma^2\,\dot z_s\,(1-\log |a_s|)\over a_s}\, ds $$ } 
$$ =\, {\sigma\over a_s} \Big[  \log |a_s| \, dM^z_s + (1-\log |a_s|)\, {\dot z_s\over a_s}\, dM^a_s\Big] + {\sigma^2\,\dot z_s\over a_s}\, ds - {\sigma^2\,\dot z_s\over 2\,a_s^3}\,(2\log |a_s| -3)\, ds \, ; $$ 
\vbox{
$$ d \bigg[ \Big({\dot z_s\over a_s}-\ell\Big) {\log |a_s|}\bigg] =\, \sigma\, d\bar M_s + \sigma^2 \Big({\dot z_s\over a_s}-\ell\Big) ds - {\sigma^2\,\ell^2\over 2\,a_s^2}\, ds  - {\sigma^2\,\dot z_s\over 2\,a_s^3}\,(2\log |a_s| -3)\, ds  \, , $$ } 
with 
$$ d\bar M_s\, :=\, a_s\1\Big[ \log |a_s| \, dM^z_s + \Big(\Big({\dot z_s\over a_s}-\ell\Big) - {\dot z_s\over a_s}\,\log |a_s|\Big) dM^a_s\Big] $$ 
and then 
$$\langle d\bar M_s\rangle = \Bigg( \bigg[{\log |a_s|\over a_s}\bigg]^2 - \bigg[{\dot z_s\,\log |a_s|\over a_s^2}\bigg]^2 + \bigg[{\dot z_s\over a_s}-\ell\bigg]\times\bigg[ 1-a_s\2 + {2\,\dot z_s\over a_s^3}\bigg] \Bigg) ds\, . $$ 
By Lemma \ref{lem.a} again, and by the almost sure integrability of $(\dot z_s/a_s-\ell)$, this implies the wanted almost sure convergence of $(\dot z_s/a_s-\ell) \log|a_s|\,$. $\;\diamond$  
}\fi 
\par \smallskip 
Propositions \ref{pro.ell}, \ref{pro.az''} and Notation (11) imply at once the following. 
\bcor  \label{cor.az} \  We have almost surely : \quad ${\ds A_s = \sqrt{1-\ell^2} + o(|a_s|^{-7/9}) . }$
\ecor 


\subsection{Study of the four-dimensional sub-diffusion $(x_s,\dot x_s,a_s,b_s)$} \label{SubDiff4} \indf 
  Recall from Section \ref{sec.RedDim} that we have  \quad 
${\ds d b_s\, =\, {\ts{3\,\sigma^2\over 2}}\, b_s\,ds + \sigma \, dM^b_s  }\,$, \quad   and  
$$ d\dot x_s\, =\, {\omega\over\sqrt{2}}\, \Big( e^{-\sqrt{2}\,\omega\,x_s}\, {b_s\over a_s} -2\Big)\, e^{-\sqrt{2}\,\omega\,x_s}\,a_s\,b_s\, ds + {\ts{3\,\sigma^2\over 2}}\,\dot x_s\,ds + \sigma\, dM^x_s\, , $$ 
the $\,\R^3$-valued martingale  $\breve M_s:= (M^a_s,M^b_s, M^x_s)$ having quadratic covariation matrix : 
$$ ((\breve K^{ij}_s)) = \pmatrix{ a_s^2 -1 & a_s\,b_s - 2\,e^{\sqrt{2}\,\omega\,x_s} & a_s\,\dot x_s\cr 
a_s\,b_s - 2\,e^{\sqrt{2}\,\omega\,x_s} & b_s^2 - 2\,e^{2\sqrt{2}\,\omega\,x_s} & b_s\,\dot x_s \cr 
a_s\,\dot x_s & b_s\,\dot x_s & \dot x_s^2 + 1 \cr } . $$ 

   The process $(b_s)$ alone is easily handled, analogously to Lemma \if{s \ref{lem.a} and }\fi   \ref{lem.z}. 
\blem \label{lem.b} \  \if{$(i)$ \}\fi   There exists a real standard real Brownian motion $(w''_s)$, and an almost surely converging process $(\eta''_s)$, such that we have : 
$$ |b_s| = \exp\Big[ \sigma^2 s + \sigma\, w''_s + \eta''_s\Big] \quad \hbox{ for any proper time }\, s \, . $$ 
\elem  
\ub{Proof} \    \if{$(i)$ \}\fi    We already noticed in Remark \ref {rem.pseudon} that the unit pseudo-norm relation $(00')$ forbids any vanishing of $\,b_s\,$. Thus there exists a real standard real Brownian motion $(w''_s)$ such that for any proper time $\,s\,$ we have :
$$ d \log |b_s| = \sigma^2 ds + \sigma^2\,e^{2\sqrt{2}\,\omega\,x_s}\,b_s\2\, ds + \sigma \, \sqrt{1 - 2\,e^{2\sqrt{2}\,\omega\,x_s}\,b_s\2}\, dw''_s \, ,  $$ 
so that there exists another real Brownian motion $\,\tilde w''\,$ such that : 
\if{ 
$$ \log |b_s| = \log |b_0| + \sigma^2 s +  \sigma^2 \int_0^s  e^{2\sqrt{2}\,\omega\,x_u}\,b_u\2\, du + \sigma\, \tilde w''\Bigg[\int_0^s \Big( 1- 2\,e^{2\sqrt{2}\,\omega\,x_u}\,b_u\2\Big) du \Bigg] \ge  \sigma^2 s + o(s) \, , $$ 
showing that $\,|b_s|\to +\ii\,$ almost surely, and even that \  ${\ds \liminf_{s\to\ii} {\log |b_s|\over \sigma^2 s} \ge 1\,}$. \parn 
This same last equation shows also that \ ${\ds \limsup_{s\to\ii} {\log |b_s|\over \sigma^2 s} \le 3/2\,}$.  \par 
   More precisely, we have : 
}\fi 
$$ \log |b_s| = \log |b_0| + \sigma^2 s +  \sigma^2 \int_0^s  e^{2\sqrt{2}\,\omega\,x_u}\,{du\over b_u^2} + \sigma\, w''_s + \sqrt{2}\,\sigma\, \tilde w''\Bigg(\!\int_0^s\! \Bigg[ { e^{2\sqrt{2}\,\omega\,x_u} \over b_u^2+ \sqrt{b_u^2- 2\,e^{2\sqrt{2}\,\omega\,x_u}} }\Bigg] du \!\Bigg) . $$ 
Now using Remark \ref {rem.pseudon} and Lemma \ref{lem.a} we get : 
$$ \int_0^\ii  e^{2\sqrt{2}\,\omega\,x_u}\,{du\over b_u^2} =  \int_0^\ii  \Big[e^{\sqrt{2}\,\omega\,x_u}\,{a_u\over b_u}\Big]^2 a_u\2\,du \,<\ii\, , $$ 
which shows that \  $ \eta''_s := \log |b_s| - \sigma^2 s - \sigma\, w''_s \,$ almost surely converges as $\,s\to\ii\,$. 
$\;\diamond $ \par\medskip 


   Then we get easily a new asymptotic variable of the relativistic diffusion. 
\blem \label{lem.ba} \  The process $\,\log(b_s/a_s)$ converges almost surely as $\,s\to\ii\,$. 
\elem  
\ub{Proof} \quad  Recalling from Remark \ref {rem.pseudon} that $\,b_s/a_s >0$, we have for any proper time $\,s\ge 0$ : 
$$ d \log \Big[{b_s\over a_s}\Big] =\, \sigma^2\,e^{2\sqrt{2}\,\omega\,x_s}\,b_s\2\, ds - \5\,\sigma^2\, a_s\2\, ds + \sigma \, ( b_s\1 d M^b_s - a_s\1 d M^a_s) \, , $$ 
or equivalently, for some real Brownian motion $\,W $ : 
$$ \log \Big[{b_s\over a_s}\Big] - \log \Big[{b_0\over a_0}\Big] =\, \sigma^2\! \int_0^s e^{2\sqrt{2}\,\omega\,x_u}\,{du\over b_u^2} - \sigma^2\! \int_0^s {du\over 2\, a_u^2} + \sigma \,W\Big( \int_0^s\! \Big[ 4\,{e^{\sqrt{2}\,\omega\,x_u}\over a_u\,b_u} - {e^{2\sqrt{2}\,\omega\,x_u}\over b_u^2} - {1\over a_u^2} \Big] du \Big) . $$ 
Now, as we already noticed in the proof of Lemma \ref{lem.b} for ${\ds \, \int_0^s  e^{2\sqrt{2}\,\omega\,x_u}\,{du\over b_u^2}\,}$, by Remark \ref {rem.pseudon} and Lemma \ref{lem.a} we get : 
$$ \int_0^\ii  {e^{\sqrt{2}\,\omega\,x_u}\over a_u\,b_u}\, du =  \int_0^\ii  \Big[e^{\sqrt{2}\,\omega\,x_u}\,{a_u\over b_u}\Big] a_u\2\,du \,<\ii\, , $$ 
and then the almost sure convergence of $\,\log ({b_s/a_s})\,$ as $\,s\to\ii\,$. $\;\diamond$  
\par \medskip 
\if{ 
$$ d \Big[{b_s\over a_s}\Big]\, =\, {d b_s\over a_s} - {b_s\, d a_s\over a_s^2} + {b_s\, \langle d a_s\rangle \over a_s^3} - {\langle d a_s, d b_s\rangle \over a_s^2} $$ 
}\fi 
By Remark \ref{rem.pseudon} again, we deduce at once the following.  
\bcor \label{cor.x} \  The process $(x_s)$ is almost surely bounded. Setting \  ${\ds \rr := \lim_{s\to\ii}\, {b_s\over a_s}\,}$, we have precisely : 
$$ 0 < \Big(1-{1\over \sqrt{2}}\Big)\,\rr\, \le\, \liminf_{s\to\ii}\, e^{\sqrt{2}\,\omega\,x_s}\le \limsup_{s\to\ii}\, e^{\sqrt{2}\,\omega\,x_s}\, \le\, \Big(1+{1\over \sqrt{2}}\Big)\,\rr \, <\ii\, . $$
\ecor 

   The following statement, analogous to Proposition \ref{pro.a/z}, ensures that the range of possible limits $\,\rr\,$ in Corollary \ref{cor.x} (and Lemma \ref{lem.ba}), is the whole $\, ]0,\ii [\,$. This provides another continuum of non-trivial bounded harmonic functions for the relativistic operator $\LL\,$. 
\bpro \label{pro.a/b} \  For any real $\,\rr_0\,$ such that $\,0< \rr_0 < \ii\,$, and for any $\,\e>0\,$, we have \quad  ${\ds\,\P \Big[\rr_0-\e < \rr = \lim_{s\to\ii}\,{b_s\over a_s}\, < \rr_0+\e \Big] > 1-\e\,}$, provided \  ${b_0/a_0}\,$ is close enough from $\,\rr_0\,$ and $\,|a_0|\,$ is large enough.
\epro 
\ub{Proof} \quad  Fix $\,A> 1\,$, $\,|a_0| > A^2$, and  $\,y_0\,$ such that \  $Y_0\,$ is close to  $\,\rr_0\,$, and use the expression displayed for $\,\log(b_s/a_s)$ in the proof of Lemma \ref{lem.ba}, Remark \ref{rem.pseudon}, and Lemma \ref{lem.z}, to get on an event of probability $>1-1/A$ :  
$$ \Big|\log \rr - \log \Big[{b_0\over a_0}\Big]\Big| \le\, \sigma^2(1+{\ts{1\over \sqrt{2}}})^2 \int_0^\ii  {du\over a_u^2} + \sigma \,\max\bigg\{ |W_u|\,\Big|\,0\le u\le 4 (1+{\ts{1\over \sqrt{2}}})\int_0^\ii  {du\over a_u^2} \bigg\}  $$ 
$$ \le \, 3\, A\2 + \sigma \,\max\{ |W_u|\,|\,0\le u\le 7\,\sigma\2 A\2\} \, , $$ 
so that \quad   ${\ds \P\Big(\, |\log \rr - \log \rr_0\,| \le\, \Big|\log\Big[{b_0\over a_0\, \rr_0}\Big]\Big|+ 3\, A\2 + A^{-1/2} \Big) > \, 1-1/A- e^{-A/14}\,}$. $\;\diamond$ 
\par \medskip 

   Recall that we set (in Section \ref{SubDiff2}) : 
$$ (11) \qquad  A_s = \sqrt{\Big[{\dot x_s\over a_s}\Big]^2 + \5\, \Big[ 2 - e^{-\sqrt{2}\,\omega\,x_s}\, {b_s\over a_s} \Big]^2}\,  =\, \sqrt{1 - {1+\dot z_s^2\over a_s^2}}\, =  \bigg|{\sh \lambda_s\over a_s}\bigg| \, . $$ 
Let us set also : 
$$ (12) \qquad  \dot x_s = a_s \, A_s\, \cos\gamma_s \; ; \quad  2 - e^{-\sqrt{2}\,\omega\,x_s}\, {b_s\over a_s}\, =\, \sqrt{2}\,A_s\, \sin\gamma_s\, . $$ 

\brem \label{rem.gamma} \quad   {\rm  By Corollary \ref{cor.az}, we have almost surely : 
$$ e^{-\sqrt{2}\,\omega\,x_s}\,{b_s\over a_s} \, = {2} - {\sqrt{2 (1-\ell^2)}}\, \sin\gamma_s + o(|a_s|^{-7/9})\; \hbox{ and } \; \dot x_s = \sqrt{1-\ell^2}\,\,a_s\,\cos\gamma_s + o(|a_s|^{2/9})\, . $$ 
}\erem 

\bpro \label{pro.alpha} \  There exists a standard real Brownian motion $\,W\,$ such that (provided $\,A_0>0$) we have almost surely, for any $\,s\ge 0$ : 
$$ \gamma_s = \gamma_0 + \omega \int_0^s e^{-\sqrt{2}\,\omega\,x_u}\, b_u\, du + \sigma\,  W\Big[\int_0^s {du\over A_u^2\,a_u^2} \Big] \, .  $$ 
\epro 
\ub{Proof} \quad   Let us differentiate the relation  
$$ \cotg\gamma_s = \sqrt{2}\,\,{\dot x_s\over a_s}\, \Big[ 2 - e^{-\sqrt{2}\,\omega\,x_s}\, {b_s\over a_s} \Big]\1\, , $$ 
which makes sense almost surely at least for a dense subset of positive $\,s\,$. \  We get : 
$$ (1+\cotg^2\gamma_s) (- d\gamma_s + \cotg\gamma_s\, \langle d\gamma_s\rangle)\, =\, \sqrt{2}\,\Big[ 2 - e^{-\sqrt{2}\,\omega\,x_s}\, {b_s\over a_s} \Big]\1 d \Big[{\dot x_s\over a_s}\Big]  \qquad $$
$$ \qquad  +\, \sqrt{2}\,\,{ {\dot x_s\over a_s}\, d\Big[ e^{-\sqrt{2}\,\omega\,x_s}\, {b_s\over a_s} \Big] + \Big\langle d\Big[{\dot x_s\over a_s}\Big] , d\Big[ e^{-\sqrt{2}\,\omega\,x_s}\, {b_s\over a_s} \Big] \Big\rangle \over \Big[ 2 - e^{-\sqrt{2}\,\omega\,x_s}\, {b_s\over a_s} \Big]^2} + \sqrt{2}\,\, {\dot x_s\over a_s}\,{ \Big\langle d\Big[ e^{-\sqrt{2}\,\omega\,x_s}\, {b_s\over a_s} \Big] \Big\rangle\over \Big[ 2 - e^{\sqrt{2}\, \omega\,x_s}\, {b_s\over a_s} \Big]^{3}}\, ,   $$ 
with 
$$ d \Big[{\dot x_s\over a_s}\Big]\, =\, {d \dot x_s\over a_s} - {\dot x_s\, d a_s\over a_s^2} - {\langle d a_s, d \dot x_s\rangle \over a_s^2} + {\dot x_s\, \langle d a_s\rangle \over a_s^3}  $$ 
$$ =\, {\omega\over\sqrt{2}}\, \Big[ e^{-\sqrt{2}\,\omega\,x_s}\, {b_s\over a_s} -2\Big]\, e^{-\sqrt{2}\,\omega\,x_s}\,b_s\, ds - \sigma^2\, {\dot x_s\over a_s^3}\, ds + {\sigma\over a_s} \,\Big[ dM^x_s - {\dot x_s\over a_s}\, dM^a_s\Big] \, , $$ 
and 
$$ d\Big[ e^{-\sqrt{2}\,\omega\,x_s}\, {b_s\over a_s} \Big] = \Big[ {2\sigma^2\over a_s^2} - e^{-\sqrt{2}\,\omega\,x_s} {\sigma^2 b_s\over a_s^3} -\sqrt{2}\,\omega\,e^{-\sqrt{2}\,\omega\,x_s}  {b_s\over a_s}\,\dot x_s \Big] ds  + {\sigma\over a_s}\,e^{-\sqrt{2}\,\omega\,x_s} \Big[ dM^b_s - {b_s\over a_s}\, dM^a_s\Big] . $$ 
Hence we get \par
\vbox{
$$ \sqrt{2}\,A_s^2\,(d\gamma_s - \cotg\gamma_s\, \langle d\gamma_s\rangle) = \Big[ e^{-\sqrt{2}\, \omega\,x_s} {b_s\over a_s} -2 \Big] \bigg[ \Big[ e^{-\sqrt{2}\,\omega\,x_s}\, {b_s\over a_s} -2\Big] e^{-\sqrt{2}\,\omega\,x_s} {\omega\,b_s\over\sqrt{2}} - \sigma^2\, {\dot x_s\over a_s^3}\bigg]  ds $$
$$ + \Big[ e^{ -\sqrt{2}\, \omega\,x_s} {b_s\over a_s} -2 \Big] {\sigma\over a_s} \,\Big[ dM^x_s - {\dot x_s\over a_s}\, dM^a_s\Big] - \bigg[ {2\sigma^2 \dot x_s\over a_s^3} - e^{-\sqrt{2}\,\omega\,x_s} {\sigma^2 b_s \dot x_s\over a_s^4} -\sqrt{2}\,\omega\,e^{-\sqrt{2}\,\omega\,x_s}  {b_s\dot x_s^2\over a_s^2} \bigg] ds $$
$$  -\, {\sigma\,\dot x_s\over a_s^2}\,e^{-\sqrt{2}\,\omega\,x_s} \Big[ dM^b_s - {b_s\over a_s}\, dM^a_s\Big] - {\sigma^2\over a_s^2}\,e^{-\sqrt{2}\,\omega\,x_s}  \Big[ 2\, e^{\sqrt{2}\,\omega\,x_s}\, {\dot x_s\over a_s} - {b_s\,\dot x_s\over a_s^2} \Big]\, ds $$ 
$$ -\, \Big[ 2 - e^{-\sqrt{2}\,\omega\,x_s}\, {b_s\over a_s} \Big]\1\,{\sigma^2\,\dot x_s\over a_s^3}\,  e^{-2\sqrt{2}\,\omega\,x_s}\,\bigg[ 4\, e^{\sqrt{2}\,\omega\,x_s}\, {b_s\over a_s} - {b_s^2\over a_s^2}  - 2\,e^{2\sqrt{2}\,\omega\,x_s} \bigg] ds $$ } \par
\vbox{
$$ =\, {\sigma\over a_s}\bigg[\Big[ e^{ -\sqrt{2}\, \omega\,x_s} {b_s\over a_s} -2 \Big]  \,dM^x_s - {\dot x_s\over a_s}\,e^{-\sqrt{2}\,\omega\,x_s}\, dM^b_s + 2 \,{\dot x_s\over a_s}  \, dM^a_s \bigg]  + \sqrt{2}\,\omega\,e^{-\sqrt{2}\,\omega\,x_s}\,b_s \,A_s^2\, ds  $$ 
$$ + \Big[ e^{-\sqrt{2}\,\omega\,x_s} {b_s\over a_s} -2 \Big] {\sigma^2 \dot x_s\over a_s^3} ds + \Big[ e^{-\sqrt{2}\,\omega\,x_s} {b_s\over a_s}-2 \Big]\1\, {\sigma^2\,\dot x_s\over a_s^3}\, \Big[ 4\, e^{-\sqrt{2}\,\omega\,x_s}\, {b_s\over a_s} - e^{-2\sqrt{2}\,\omega\,x_s}\,{b_s^2\over a_s^2} - 2\Big]\, ds \, . $$ }
This implies 
$$ 2\,A_s^4\,\langle d\gamma_s\rangle\, =\,  {\sigma^2\over a_s^2}\,\bigg\langle\Big[ e^{ -\sqrt{2}\, \omega\,x_s} {b_s\over a_s} -2 \Big]  \,dM^x_s - {\dot x_s\over a_s}\,e^{-\sqrt{2}\,\omega\,x_s}\, dM^b_s + 2 \,{\dot x_s\over a_s}  \, dM^a_s \bigg\rangle = 2\sigma^2\,a_s\2\, A_s^2\, , $$ 
whence for some standard real Brownian motion $\,\tilde W$ : \par \smallskip 
\vbox{
$$ d\gamma_s\,=\, {\sigma\over A_s\,a_s}\, d\tilde W_s + \omega\,e^{-\sqrt{2}\,\omega\,x_s}\,b_s \, ds + \Big[ e^{-\sqrt{2}\,\omega\,x_s} {b_s\over a_s} -2 \Big]\, {\sigma^2\, \dot x_s\over \sqrt{2}\,A_s^2\,a_s^3}\, ds $$ 
$$ - \Big[ e^{-\sqrt{2}\,\omega\,x_s} {b_s\over a_s}-2 \Big]\1\, {\sigma^2\, \dot x_s\over  \sqrt{2}\, A_s^2\,a_s^3} \bigg[ 2 - \Big[ 4\, e^{-\sqrt{2}\,\omega\,x_s}\, {b_s\over a_s} - e^{-2\sqrt{2}\,\omega\,x_s}\,{b_s^2\over a_s^2}  - 2\Big] \bigg] ds $$ } \par
\vspace{-3mm} 
$$ =\,  {\sigma\over A_s\,a_s}\, d\tilde W_s + \omega\,e^{-\sqrt{2}\,\omega\,x_s}\,b_s \, ds \, . $$ 
\par \smallskip \noindent 
Recall from Remark \ref{rem.lambda}  that (provided $\,A_0>0\,$) $\, A_s\,a_s\,$ almost surely never vanishes. This shows that the above formula almost surely holds for any $\,s\ge 0\,$, proving the result.  $\;\diamond $
\par\medskip 

 By Proposition \ref{pro.a-z} and by the definition $(11)$ of $\,A_s\,$,  we have almost surely : 
$$ A_u^2\,a_u^2\, =\,\sh^2\lambda_u \, =\,\sh^2[\sigma^2 u + \sigma\, \tilde w_u + \tilde\eta_u]\, , $$
which proves the almost sure convergence of \ ${\ds \int^\ii {du\over A_u^2\,a_u^2}}\,$. Hence we deduce the following. 
\bcor \label{cor.gamma} \  There exists a converging process $\,\check\eta\,$ such that we have almost surely, for any $\,s\ge 0$ : 
$$ \gamma_s = \omega \int_0^s e^{-\sqrt{2}\,\omega\,x_u}\, b_u\, du + \check\eta_s \, .  $$ 
\ecor 

Using Corollary \ref{cor.gamma}, Corollary \ref{cor.x}, Lemma \ref{lem.b} and Remark \ref{rem.gamma}, we deduce also easily the following.  
\bcor \label{cor.gamma'} \  We have almost surely, for large $\,s$ : \  
${\ds  \gamma_s \, =\, {\rm sign}(b_0)\,  e^{\sigma^2 s+ o(s^{5/9})}\, , }$ 
$$ e^{-\sqrt{2}\,\omega\,x_s} = {\ts{2\over\rr}} - {\ts{\sqrt{2 (1-\ell^2)} \over {\rm sign}(b_0) \rr}}  \sin\!\Big( e^{\sigma^2s+o(s^{5\over 9})}\Big) + o\Big(\! e^{-{2\over 3}\sigma^2 s}\!\Big)  , \; {\dot x_s\over a_s} = \sqrt{1-\ell^2}\,\cos\!\Big( e^{\sigma^2s+o(s^{5\over 9})}\Big) + o\Big(\! e^{-{2\over 3}\sigma^2 s}\!\Big)  . $$
\ecor 

   In the same vein as Proposition \ref{pro.az''}, we have the following. 
\bpro \label{pro.ab'} \   We have \ $(b_s - \rr\,a_s) = o(|a_s|^{2/9})$, almost surely as $\,s\to\ii$. 
\epro 
\ub{Proof} \quad  We have : \quad 
${\ds d \bigg[{b_s\over a_s}\bigg] =\, 2\sigma^2\,{e^{\sqrt{2}\,\omega\,x_s}\over a_s^2}\, ds - \sigma^2\, {b_s\over a_s^3}\, ds + {\sigma\over a_s} \, \Big[d M^b_s - {b_s\over a_s}\, d M^a_s\Big] \, , }$ 
$$  \Big\langle{ dM^b_s- {b_s\over a_s}\,dM^a_s}\Big\rangle = \Big( 4\, e^{\sqrt{2}\,\omega\,x_s}
{b_s\over a_s} - 2\, e^{2\sqrt{2}\,\omega\,x_s} - {b_s^2\over a_s^2}
\Big) ds\; , \quad \Big\langle{dM^a_s\over a_s}\Big\rangle = (1-a_s\2) ds \, ,  $$ 
and
$$ \Big\langle{ {dM^a_s}\, ,\,dM^b_s- {b_s\over a_s}\,dM^a_s}\Big\rangle = \Big({b_s\over a_s} - 2\, e^{\sqrt{2}\,\omega\,x_s}\Big) ds \, . $$
Hence there exist two independent standard real Brownian motions $\,w,w'\,$ such that : 
$$ dM^b_s- {b_s\over a_s}\,dM^a_s = {b_s - 2\, a_s\,e^{\sqrt{2}\,\omega\,x_s}\over a_s^2\,\sqrt{1-a_s\2}}\,dw_s + \sqrt{2\, e^{2\sqrt{2}\,\omega\,x_s}(\dot z_s^2+\dot x_s^2)\over a_s^2-1}\,\, dw'_s\, . $$ 
and as in the proof of Proposition \ref{pro.az''} : 
$$  {dM^a_s\over a_s} = \sqrt{1-a_s\2}\, dw_s \, ,\hbox{ and }\quad {a_s\over a_{s+u}} \,=\, e^{-\sigma^2 u - \sigma (w_{s+u} -w_s) +  \O(a_s\1)} . $$  

   Therefore, we get almost surely : \par
\vbox{
$$ \rr - {b_s\over a_s}\, =\, 2\sigma^2\!\int_s^\ii {e^{\sqrt{2}\,\omega\,x_u}\over a_u^2}\, du - \sigma^2\!  \int_s^\ii {b_u\over a_u^3}\, du + \sigma\!\int_s^\ii  {b_u - 2\, a_u\,e^{\sqrt{2}\,\omega\,x_u}\over a_u^3\,\sqrt{1-a_u\2}}\,dw_u \quad $$ 
$$ \qquad +\, \sigma\!\int_s^\ii {e^{\sqrt{2}\,\omega\,x_u}\over a_u} \sqrt{2\,(\dot z_u^2+\dot x_u^2)\over a_u^2-1}\,\, dw'_u  \, . $$ } \parn  

Now, as in the proof of Proposition \ref{pro.az''}, we have 
$$ \int_s^\ii {e^{\sqrt{2}\,\omega\,x_u}\over a_u^2}\, du + \int_s^\ii {b_u\over a_u^3}\, du  = o\Big(|a_s|^{-7/4}\Big) \, , \hbox{ and } \int_s^\ii  {b_u - 2\, a_u\,e^{\sqrt{2}\,\omega\,x_u}\over a_u^3\,\sqrt{1-a_u\2}}\,dw_u = o\Big(|a_s|^{-4/3}\Big) , $$ 
and 
$$\int_s^\ii {e^{\sqrt{2}\,\omega\,x_u}\over a_u} \sqrt{2\,(\dot z_u^2+\dot x_u^2)\over a_u^2-1}\,\, dw'_u = o\Big(|a_s|^{-7/9}\Big) \, . $$

   So far, we have shown that almost surely, as $\,s\to\ii$ : 
$$ {b_s} -\rr\,a_s\,  = \, o\Big(|a_s|^{-1/3}\Big) - \sigma\, a_s\! \int_s^\ii {e^{\sqrt{2}\,\omega \,x_u} \over a_u} \sqrt{2\,(\dot z_u^2+\dot x_u^2)\over a_u^2-1}\,\, dw'_u\, =\,  o\Big(|a_s|^{2/9}\Big) . 
\;\;\diamond $$ 

\brem \label{rem.ab'} \  {\rm  Almost the same proof as above shows that we have in fact \ $(b_s - \rr\,a_s) = o(|a_s|^{\e})$ almost surely, and similarly, \ $(\dot z_s - \ell\,a_s) = o(|a_s|^{\e})$  almost surely, for any $\,\e>0$. But Propositions \ref{pro.az''}, \ref{pro.ab'} (as Corollaries \ref{cor.compz}, \ref{cor.az}, \ref{cor.gamma'}) as stated are sufficient for our purpose. \par \smallskip

   It is possible to complete Proposition \ref{pro.ab'} in a way somewhat similar to Proposition \ref{pro.az''} : starting from the last formula above (in the proof of Proposition \ref{pro.ab'}), and proceeding as in the proof of Proposition \ref{pro.az''}, we can see that the process $(b_s - \rr\,a_s)$ asymptotically behaves as 
$$ \sigma\, a_s\! \int_s^\ii {e^{\sqrt{2}\,\omega \,x_u} \over a_u} \sqrt{2\bigg[{\dot z_u^2+\dot x_u^2\over a_u^2}\bigg]}\,\, dw''_u\,, $$  
for a standard real Brownian motion $\,w''$, independent from the filtration of $(a_s, b_s, \dot x_s, \dot z_s)$, and then, using the pseudo-norm relation $(00')$ together with Proposition \ref{pro.ab'}, as 
$$ \sigma a_s\! \int_s^\ii\! \sqrt{4 \rr\, e^{\sqrt{2}\,\omega x_u} - 2\, e^{2\sqrt{2}\,\omega x_u} -\rr^2}\, {dw''_u\over a_u} = W''\bigg[ \!\int_0^\ii\! \Big[4 \rr e^{\sqrt{2}\,\omega x_{s+u}} - 2 e^{2\sqrt{2}\,\omega x_{s+u}} -\rr^2\Big] {\sigma^2a_s^2\over a_{s+u}^2}\, du \bigg] . $$  
Note that to complete the Brownian representation of the martingale $(\breve M_s)$ in the proof of Proposition \ref{pro.ab'}, we can find a Brownian motion $\,w'''$ in the filtration of $(a_s, b_s, \dot x_s, \dot z_s)$, independent from $\,w,w'$, such that 
$$ dM^x_s = {a_s\,\dot x_s\over \sqrt{a_s^2-1}}\,dw_s - {(b_s - 2\, a_s\,e^{\sqrt{2}\,\omega\,x_s}) \, \dot x_s\over \sqrt{a_s^2-1}\,\sqrt{2\, e^{2\sqrt{2}\,\omega\,x_s}(\dot z_s^2+\dot x_s^2)}}\,\, dw'_s + {\dot z_s\over \sqrt{\dot z_s^2+\dot x_s^2}}\, dw'''_s\, . $$ 
Now \  
$$ \bigg( \rr\, , \int_0^\ii\! e^{\sqrt{2}\,\omega (x_{s+u}-x_{s})} \, {a_s^2\over a_{s+u}^2}\, du \,,\; \int_0^\ii\! e^{2\sqrt{2}\,\omega (x_{s+u}-x_{s})} \, {a_s^2\over a_{s+u}^2}\, du\,,\; \int_0^\ii\! {a_s^2\over a_{s+u}^2}\, du\,,\; e^{\sqrt{2}\,\omega\, x_{s}}\bigg) $$ 
converges in law to \  
$$ \bigg( \rr\, , \int_0^\ii\! e^{\sqrt{2}\,\omega\, x_{u}-2\sigma^2u-2\sigma w_u} \, du\,, \int_0^\ii\! e^{2\sqrt{2}\,\omega\, x_{u}-2\sigma^2u-2\sigma w_u} \, du\,, \int_0^\ii\! e^{-2\sigma^2u-2\sigma w_u} \, du\,,\; V\bigg) , $$ 
where \quad $V\1:= {\ts{2\over\rr}} - {\ts{\sqrt{2 (1-\ell^2)} \over \rr}} \sin U\,$, $\,U\,$ denoting an independent variable, uniform on the circle.  We get thus the convergence in law, but not in probability, of $(b_s - \rr\,a_s)$, to the smooth : 
$$ W''\bigg[ \sigma^2\!\int_0^\ii\! \Big[4 \rr\,V e^{\sqrt{2}\,\omega\, x_{u}} - 2\,V^2 e^{2\sqrt{2}\,\omega \, x_{u}} -\rr^2\Big] e^{-2\sigma^2u-2\sigma w_u}\, du \bigg] .
$$ 
}\erem

\subsection{Irreducibility} \label{sec.irred}  

\bpro \label{pro.irred} \  $(i)$ \  The relativistic diffusion is irreducible : from any starting point, it hits any non-empty open subset of the phase-space $\,\EE\,\moins\,\EE_0\,$ with a strictly positive probability.  \parn 
$(ii)$ \  For any starting point (in $\,\EE$), the law of the asymptotic variable $(\ell,\rr)$ charges any non-empty open subset of the range $\Big( ]-1,0[\,\cup\, ]0,1[\Big)\times ]0,\ii[\,$. 
\epro
\ub{Proof} \ $(i)$ \    We know from Proposition \ref{pro.geod} (in Section \ref{sec.geod}) that there are piece-wise geodesic timelike continuous paths, and then trajectories in the support of the relativistic diffusion $(\xi_\cdot,\dot\xi_\cdot)$, moving at will the coordinates $(t,x,y,z)$. \par 

   Owing to the the quadratic covariation (rank 3) matrix of the $\,\R^3$-valued martingale   $(M^a_s,M^b_s, M^z_s)$ (recall Section \ref{sec.RedDim}) and to the unit pseudo-norm relation $(00)$, we can find three independent standard real Brownian motions $(w^1,w^2,w^3)$ such that : 
$$ d \dot z_s\, =\, {\ts{3\,\sigma^2\over 2}}\, \dot z_s\,ds + \sigma\, \sqrt{\dot z_s^2+1}\,\, dw^1_s \; ; $$ 
$$ d a_s\, =\, {\ts{3\,\sigma^2\over 2}}\, a_s\,ds + \sigma\,{a_s\,\dot z_s\over \sqrt{\dot z_s^2+1}}\, dw^1_s + \sigma\,\sqrt{{a_s^2-\dot z_s^2-1\over \dot z_s^2+1}}\,\, dw^2_s\; ; $$ 
$$ d b_s = {\ts{3\,\sigma^2\over 2}}\, b_s\,ds + \sigma\,{b_s\,\dot z_s\over \sqrt{\dot z_s^2+1}}\, dw^1_s + \sigma\, {a_s\,b_s- 2\,e^{\sqrt{2}\,\omega\,x_u}(\dot z_s^2+1) \over \sqrt{(\dot z_s^2+1)(a_s^2-\dot z_s^2-1)}}\, dw^2_s + \sigma\,{ \sqrt{2}\, e^{\sqrt{2}\,\omega\,x_u}\,\dot x_s \over \sqrt{a_s^2-\dot z_s^2-1} }\, dw^3_s\, . $$ 

   Let us use the support theorem of Stroock and Varadhan (see for example ([I-W], Theorem VI.8.1)). We see thus from the above stochastic differential system, that the following trajectories belong to the support of $(\xi_\cdot,\dot\xi_\cdot) \equiv (t_\cdot,x_\cdot,y_\cdot,z_\cdot, \dot z_\cdot,a_\cdot,b_\cdot, \dot x_\cdot)$ : \parn 
- trajectories moving at will the coordinate $\,\dot z\,$, without changing the coordinates $(t,x,y,z)$ ; \parn 
- trajectories moving at will the coordinate $\,a\,$, without changing the coordinates \parn $(t,x,y,z, \dot z)$ ; \parn 
- trajectories moving at will the coordinate $\,b\,$, provided $\,\dot z\not =0\,$,  without changing the coordinates $(t,x,y,z, \dot x,a)$. \par  \smallskip 

  So far, it has become clear that it is possible, within the support of the relativistic diffusion, to move any point of the phase space $\,\EE$ having given first coordinates $(t,x,y,z,\dot z,a)\in\R^6$, onto some point of the phase space $\,\EE$ having prescribed first coordinates $(t',x',y',z', \dot z',a')\in\R^6$. \par 
  It remains only to consider the last two coordinates $(b,\dot x)$. They are of course constrained by the unit pseudo-norm relation $(00')$, which tells precisely that they run some ellipse of this plane of coordinates, which is centred on the axis $\{\dot x=0\}$. The last type of trajectory mentioned above  allows now to move $(b,\dot x)$ arbitrarily on the upper half and on the lower half this ellipse, without changing the other coordinates, within the support of the relativistic diffusion. \par \smallskip 
  Finally, we must make clear that we can cross the axis $\{\dot x=0\}$, within the support of the relativistic diffusion, at the cost of an arbitrarily small move of all coordinates. Now, owing to the the quadratic covariation (rank 3) matrix of the $\,\R^3$-valued martingale   $(M^a_s,M^b_s, M^x_s)$ (recall Section \ref{sec.RedDim}) and to the unit pseudo-norm relation $(00)$, we can find three independent standard real Brownian motions $(\bar w^1,\bar w^2,\bar w^3)$ such that : 
$$ d \dot x_s\, =\, (\omega/\sqrt{2}\,)\, e^{-2\sqrt{2}\,\omega\,x_s}\, b_s^2\, ds - \sqrt{2}\,\omega\, e^{-\sqrt{2}\,\omega\,x_s}\,a_s\,b_s\, ds + {\ts{3\,\sigma^2\over 2}}\, \dot x_s\,ds + \sigma\, \sqrt{\dot x_s^2+1}\,\, d\bar w^1_s \; ; $$ 
$$ d a_s\, =\, {\ts{3\,\sigma^2\over 2}}\, a_s\,ds + \sigma\,{a_s\,\dot x_s\over \sqrt{\dot x_s^2+1}}\, d\bar w^1_s + \sigma\,\sqrt{{a_s^2-\dot x_s^2-1\over \dot x_s^2+1}}\,\, d\bar w^2_s\; ; $$ 
$$ d b_s = {\ts{3\,\sigma^2\over 2}}\, b_s\,ds + \sigma\,{b_s\,\dot x_s\over \sqrt{\dot x_s^2+1}}\, d\bar w^1_s + \sigma\, {a_s\,b_s- 2\,e^{\sqrt{2}\,\omega\,x_u}(\dot x_s^2+1) \over \sqrt{(\dot x_s^2+1)(a_s^2-\dot x_s^2-1)}}\, d\bar w^2_s + \sigma\,{ \sqrt{2}\, e^{\sqrt{2}\,\omega\,x_u}\,\dot z_s \over \sqrt{a_s^2-\dot x_s^2-1} }\, d\bar w^3_s\, . $$ 
This shows that, arrived in an arbitrarily thin $\,\delta$-neighbourhood of the axis $\{\dot x=0\}$, we can cross this axis (within the support of the relativistic diffusion, acting on the Brownian component $\bar w^1$) without changing the coordinates $(t,x,y,z)$, and perturbing the coordinates $(a,b,\dot x,\dot z)$ only by some move of order $\,\delta\,$. This ends the proof of  irreducibility.  \par \medskip 

   $(ii)$ \  This is a direct consequence of $(i)$ above and of Propositions \ref{pro.a/z} and \ref{pro.a/b} : by $(i)$, it is indeed enough to start the relativistic diffusion so that ${\dot z_0/a_0}\,$ be close to  a given $\,\ell_0\in \Big( ]-1,0[\,\cup\, ]0,1[\Big)$, ${b_0/a_0}\,$ be close to  a given $\,\rr_0>0\,$, and $\,|a_0|\,$ be large enough. $\;\diamond$ 
   
\subsection{Convergence to a lightlike geodesic} \label{sec.convlg} \indf 
  From Section \ref{sec.ngeod} (Proposition \ref{pro.geodlum} and Remark \ref{rem.geodlum}), for a lightlike geodesic we have three geometrically meaningful real parameters $\,\ell, \rr, Y\,$ such that : 
$$ {\dot z_\tau\over a} = \ell \,\in\,[0,1]\; ; \quad {b\over a}  = \rr \,\in \,]0,\ii[ \; ; \quad {\sqrt{2} \over \omega\, b}\,\,\dot x_\tau + y_\tau = Y \,\in\R \; , $$
$$ \hbox{and }\qquad  \Big[ {\rr\over 2}\, e^{-\sqrt{2}\,\omega\,x_\tau} - 1\Big]^2 + \Big[{\omega\, \rr\over 2}\, (y_\tau - Y)\Big]^2 =\, \5\,(1-\ell^2)\, , $$ 
where $\,\tau\,$ denotes an (irrelevant) affine parameter. \par \smallskip  

   Accordingly, let us exhibit now a third asymptotic random variable $Y$ for the relativistic diffusion. 
\bpro \label{pro.Y} \   The process \  ${\ds Y_s := \, y_s + {\ts{\sqrt{2}\over \omega}}\,{\dot x_s\over b_s}\,}$ converges almost surely, as $\,s\to\ii$, toward some real random variable $\,Y$. 
\  And we have \  $Y-Y_s= o(|a_s|^{-4/9})\,$ almost surely. 
\epro 
\ub{Proof} \quad  Recall from Formulas $(10)$ that we have \  $\dot y_s = e^{-\sqrt{2}\,\omega\,x_s} ( 2\,a_s- e^{-\sqrt{2}\,\omega\,x_s}b_s)$.   
We have then : 
$$ d \bigg[{\dot x_s\over b_s}\bigg]\, =\, {d \dot x_s\over b_s} - {\dot x_s\, d b_s\over b_s^2} + {\dot x_s\, \langle d b_s\rangle \over b_s^3} - {\langle d b_s, d \dot x_s\rangle \over b_s^2}  $$ 
$$ =\, {\ts{\omega\over\sqrt{2}}}\, e^{- 2\sqrt{2}\,\omega\,x_s}\, {b_s}\,ds - \sqrt{2}\,\omega\, e^{-\sqrt{2}\,\omega\,x_s}\,a_s\, ds - 2\sigma^2\, e^{2\sqrt{2}\,\omega\,x_s}\, {\dot x_s\over b_s^3}\, ds + {\sigma\over b_s} \, dM^x_s - \sigma \, {\dot x_s\over b_s^2}\, dM^b_s \, , $$ 
whence 
$$ {\ts{\omega\over\sqrt{2}}}\, dY_s\, = \, - 2\sigma^2\, e^{2\sqrt{2}\,\omega\,x_s}\, {\dot x_s\over b_s^3}\, ds + {\sigma\over b_s} \, dM^x_s - \sigma \, {\dot x_s\over b_s^2}\, dM^b_s \, , $$ 
and for some Brownian motion $\,W'$ : 
$$ {\ts{\omega\over\sqrt{2}}}\, Y_s\, = \, {\ts{\omega\over\sqrt{2}}}\, Y_0 - 2\sigma^2\! \int_0^s  e^{2\sqrt{2}\,\omega\,x_u}\, {\dot x_u\over b_u^3}\, du + \sigma\, W'\!\left[ \int_0^s \Big[ 1- e^{2\sqrt{2}\,\omega\,x_u}\, \Big|{\dot x_u\over b_u}\Big|^2\Big]\, {du\over b_u^2} \right] . $$ 
By Corollary \ref{cor.x}, Remark \ref{rem.gamma} and Proposition \ref{pro.ab'}, the two above integrals, and then $Y_s\,$, converge almost surely, and moreover, we have (for some Brownian motion $\,W$) :  
$$ Y_s - Y = Y_s - Y_\ii\, = \, {\ts{2\sqrt{2}\,\sigma^2\over \omega}} \int_s^\ii  e^{2\sqrt{2}\,\omega\,x_u}\, {\dot x_u\over b_u^3}\, du +  W\!\left[  {\ts{2\,\sigma^2\over \omega^2}}\! \int_s^\ii \Big[ 1- e^{2\sqrt{2}\,\omega\,x_u}\, \Big|{\dot x_u\over b_u}\Big|^2\Big]\, {du\over b_u^2} \right]  $$ 
$$ =\,\O(1)\! \int_s^\ii\! a_s\2\, ds\, +  W\!\left[  \O(1)\! \int_s^\ii\! a_u\2\, {du} \right] =\, o\left[ \! \int_s^\ii\! a_u\2\, {du} \right]^{4/9} =\, o\left( a_s^{-8/9}\right) \left[\!\int_0^\ii\! \Big|{a_s\over a_{s+u}}\Big|^2 {du} \right]^{4/9}  $$ 
$$ =\, o\left( a_s^{-8/9}\right) \left[\!\int_0^\ii\! e^{-2\sigma^2u-2\sigma (w_{s+u}-w_s)} {du} \right]^{4/9} =\, o\left( a_s^{-4/9}\right) \left[\!\int_0^\ii\! e^{-2\sigma^2u-2\sigma\,w_{s+u} -\sigma^2 s +\sigma\,w_s} {du} \right]^{4/9}  $$ 
$$ =\, o\left( a_s^{-4/9}\right) .\;\;\diamond $$ 

\if{ 
\bpro \label{pro.Y} \   The process \  ${\ds Y_s := \, y_s + {\ts{\sqrt{2}\over \omega\, \rr}}\,{\dot x_s\over a_s}\,}$ converges almost surely, as $\,s\to\ii$, toward some real random variable $\,Y$. 
\  And we have \  $Y-Y_s= o(|a_s|^{-4/9})\,$ almost surely. 
\epro 
\ub{Proof} \quad  Recall from Formulas $(10)$ that we have \  $\dot y_s = e^{-\sqrt{2}\,\omega\,x_s} ( 2\,a_s- e^{-\sqrt{2}\,\omega\,x_s}b_s)$, hence  : 
$$ y_s -y_0 = \int_0^s e^{-\sqrt{2}\,\omega\,x_u} ( 2\,a_u- e^{-\sqrt{2}\,\omega\,x_u}b_u) du\, . $$  
We have then : 
$$ d \bigg[{\dot x_s\over a_s}\bigg]\, =\, {d \dot x_s\over a_s} - {\dot x_s\, d a_s\over a_s^2} + {\dot x_s\, \langle d a_s\rangle \over a_s^3} - {\langle d a_s, d \dot x_s\rangle \over a_s^2}  $$ 
$$ =\, {\omega\over\sqrt{2}}\, \Big( e^{-\sqrt{2}\,\omega\,x_s}\, {b_s\over a_s} -2\Big)\, e^{-\sqrt{2}\,\omega\,x_s}\,b_s\, ds - \sigma^2\, {\dot x_s\over a_s^3}\, ds + {\sigma\over a_s} \, dM^x_s - \sigma \, {\dot x_s\over a_s^2}\, dM^a_s \, , $$ 
whence for some Brownian motion $\,W$ : 
$$  {\dot x_s\over a_s} - {\omega\over\sqrt{2}} \int_0^s \Big( e^{-\sqrt{2}\,\omega\,x_u}\, {b_u\over a_u} -2\Big)\, e^{-\sqrt{2}\,\omega\,x_u}\,b_u\, du \, =\,{\dot x_0\over a_0} - \sigma^2\!  \int_0^s {\dot x_u\over a_u^3}\, du + \sigma\, W\Big(\! \int_0^s\! \Big[ 1-\Big|{\dot x_u\over a_u}\Big|^2\Big] {du\over a_u^2} \Big) , $$ 
which converges almost surely as $\,s\to\ii$ (recall that by Equation $(00')$, $\,|{\dot x_s/ a_s}|\le 1$), and moreover, as seen in the proof of Proposition \ref{pro.a/b}, differs from its limit by $\,o(|a_s|^{-7/9})\,$ almost surely. 
Hence 
$$ \tilde Y_s := \, Y_s - \int_0^s e^{-\sqrt{2}\,\omega\,x_u} ( 2\,a_u- e^{-\sqrt{2}\,\omega\,x_u}\,b_u)\, du - {1\over\rr} \int_0^s \Big( e^{-\sqrt{2}\,\omega\,x_u}\, {b_u\over a_u} -2\Big)\, e^{-\sqrt{2}\,\omega\,x_u}\,b_u\, du  $$ 
converges almost surely as $\,s\to\ii$, with moreover \  $\tilde Y_s= \tilde Y_\ii + o(|a_s|^{-7/9})\,$. Now, we have by Remark \ref{rem.gamma} : 
$$ e^{-\sqrt{2}\,\omega\,x_u}\, {b_u\over a_u}\, =\,  {2} - {\sqrt{2 (1-\ell^2)}}\, \sin\gamma_u + o(|a_u|^{-7/9}) \, , $$ 
whence by Proposition \ref{pro.ab'} : 
$$ Y_s -\tilde Y_s =  \rr\1 \sqrt{2 (1-\ell^2)} \! \int_0^s e^{-\sqrt{2}\,\omega\,x_u}\Big(\rr -{b_u\over a_u}\Big) \Big[ a_u\, \sin\gamma_u +o(|a_u|^{2/9}) \Big] du  $$ 
$$ =\,  \rr\1 \sqrt{2 (1-\ell^2)} \! \int_0^s e^{-\sqrt{2}\,\omega\,x_u}\,(\rr\,a_u -b_u)\,\sin\gamma_u\, du\, + \int_0^s o(|a_u|^{-5/9})\, du \, . $$ 

   So far, we have got the almost sure convergence of 
$$ Y_s +  \rr\1 \sqrt{2 (1-\ell^2)} \! \int_0^s e^{-\sqrt{2}\,\omega\,x_u}\,(b_u -\rr\,a_u)\, \sin\gamma_u\, du\,  , $$ 
with (again as in the proof of Proposition \ref{pro.a/b}) a difference of order $\,o(|a_s|^{-4/9})$ with its limit. 
Now, by Proposition \ref{pro.alpha} we have almost surely : 
$$ \int_0^s e^{-\sqrt{2}\,\omega\,x_u}\,(b_u - \rr\,a_u)\,\sin\gamma_u\, du\, =  
\int_0^s \Big( 1 - \rr\,{a_u\over b_u}\Big)\,\sin\gamma_u\, \Big( \omega\1 d\gamma_u + \O(a_u\1) d\tilde W_u\Big) $$
$$ =\, o(|a_s|^{-7/9}) - \omega\1\!\int_0^s \Big( 1 - \rr\,{a_u\over b_u}\Big)\, d(\cos\gamma_u) $$
$$ = o(|a_s|^{-7/9}) - \omega\1\!\Big[\Big( 1 - \rr\,{a_u\over b_u}\Big)\cos\gamma_u\Big]_0^s - {\rr\over\omega} \int_0^s \cos\gamma_u\, d\Big[ {a_u\over b_u}\Big]   =\, o(|a_s|^{-7/9}) - {\rr\over\omega} \int_0^s \cos\gamma_u\, d\Big[ {a_u\over b_u}\Big] . $$ 
Otherwise, for some Brownian motion $\,\beta\,$  we have : 
$$ d\Big[ {a_u\over b_u}\Big] = {\sigma\over b_u}\,\sqrt{ 4\,{a_u\over b_u}\, e^{\sqrt{2}\,\omega\,x_u} - 2\, {a_u^2\over b_u^2}\, e^{2\sqrt{2}\,\omega\,x_u} - 1}\,\, d\beta_u + 2\, {e^{\sqrt{2}\,\omega\,x_u}\over b_u^2}\Big[ 1- {a_u\over b_u}\, e^{\sqrt{2}\,\omega\,x_u}\Big] du\, , $$ 
whence by Corollary \ref{cor.x} :  \  
${\ds \cos\gamma_u\, d\Big[ {a_u\over b_u}\Big] =\, \O(|b_u|\1)\, d\beta_u + \O(b_u\2)\, du \, . }$ 
 Finally, using the above and Lemma \ref{lem.b}, we get the almost sure convergence of \  ${\ds \int_0^s e^{-\sqrt{2}\,\omega\,x_u} (b_u - \rr\,a_u) \sin\gamma_u\, du\,}$, \  with a difference of order $\,o(|a_s|^{-4/9})\,$ from its limit,  and then the claim.  $\;\diamond$  
\par\medskip 
}\fi 

   Proposition \ref{pro.Y}, together with Remark \ref{rem.gamma} and Corollary \ref{cor.x}, implies easily the following. 
\bcor \label{cor.Y} \  We have almost surely, as $\,s\to\ii$ : 
$$  \Big[ {\rr\over 2}\, e^{-\sqrt{2}\,\omega\,x_s} - 1\Big]^2 + \Big[{\omega\, \rr\over 2}\, (y_s - Y)\Big]^2 \lra\, \5(1-\ell^2)\, . $$ 
\ecor 

   The following statement, analogous to Propositions \ref{pro.a/z} and \ref{pro.a/b}, ensures that the range of possible limits $\,Y\,$ in Proposition \ref{pro.Y} is the whole $\,\R\,$. This provides again another continuum of non-trivial bounded harmonic functions for the relativistic operator $\LL\,$. 
\bpro \label{pro.Y'} \  For any real $\,y\,$ and any $\,\e>0$, we have \   ${\ds\,\P [ y-\e <  Y  <  y+\e ] > 1-\e\,}$, provided \  $Y_0\,$ is close enough from $\,y\,$ and $\,|a_0|\,$ is large enough.
\epro 
\ub{Proof} \quad  Recall from Lemma \ref{lem.z} that the event \  ${\ds \AA := \Big\{ |a_s| \ge \sqrt{|a_0|}\, e^{\sigma^2 s/2} \hbox{ for any } s\ge 0\Big\} }$ has (for $|a_0|>3$) probability larger than $\, 1-|a_0|^{-1/2}$. \  The proof of  Proposition \ref{pro.Y} shows that 
$$ |Y -Y_0| = \O(1)\! \int_0^\ii {du\over a_u^2}  + \max\bigg\{ |W_s|\,\Big|\,0\le s\le \O(1)\! \int_0^\ii  {du\over a_u^2} \bigg\} = \O({\ts{1\over |a_0|}}) + W^*({\ts{1\over |a_0|}}) \; \hbox{ on } \AA\, ,  $$ 
so that \quad  ${\ds \P\Big(\, |Y -Y_0\,| \le\, 2\,|a_0|^{-1/3} \Big) > \, 1-2\,|a_0|^{-1/2}\,}$, for large enough $\,|a_0|\,$.  
\if{\  The proof of  Proposition \ref{pro.Y} shows then also that : 
$$ |Y-\tilde Y_\ii| \le C\! \int_0^\ii {du\over a_u^{4/9}}  + \max\bigg\{ |W_s|\,\Big|\,0\le s\le C\!\int_0^\ii  {du\over a_u^2} \bigg\} +  \max\bigg\{ |B_s|\,\Big|\,0\le s\le C\!\int_0^\ii  {du\over b_u^2} \bigg\}  $$ 
holds almost surely, for some constant $\,C\,$ and Brownian motions $B,W$. Whence, using Proposition \ref{pro.ab'} to handle the last term, for some other constant $\,C$ : 
$$ |Y-\tilde Y_\ii| \le C\! \int_0^\ii {du\over a_u^{4/9}}  + \max\bigg\{ |W_s|\,\Big|\,0\le s\le C\!\int_0^\ii  {du\over a_u^2} \bigg\} +  \max\bigg\{ |B_s|\,\Big|\,0\le s\le C\!\int_0^\ii  {du\over a_u^2} \bigg\} , $$ 
which holds then with high probability for some non-random constant $\,C\,$.  
Proceeding finally as for $ |\tilde Y_\ii-Y_0|$ above, we get the claim.
}\fi
$\;\diamond$ 
\par \medskip 

    Proposition \ref{pro.Y'} improves Proposition (\ref{pro.irred},$(ii)$). We deduce indeed at once the following. 
\bcor \label{cor.irred} \   For any starting point (in $\,\EE$), the law of the asymptotic variable $(\ell,\rr,Y)$ charges any non-empty open subset of the range $\Big( ]-1,0[\,\cup\, ]0,1[\Big)\times ]0,\ii[\times \R\,$. \par
   More precisely,  if the starting point of the relativistic diffusion satisfies : ${\dot z_0/a_0}\,$ close enough to  $\,\ell_0\in ]-1,1[\,$, ${b_0/a_0}\,$ close enough to $\,\rr_0>0\,$, $Y_0\,$ close enough to $\,y\in\R\,$, and $\,|a_0|\,$ large enough, then with arbitrary large probability, $(\ell,\rr,Y)$ is arbitrary close to  $(\ell_0,\rr_0,y)$.
\ecor     
    
\brem \label{rem.lum} \  {\rm A rapid look at Remark \ref{rem.geodlum} could let think that there could be a fourth asymptotic random variable for the relativistic diffusion, namely a possible almost sure limit for \  $X_s := z_s+\ell\,t_s -(\ell/\omega) \gamma_s\,$. But as a matter of fact, there is no such limit, in accordance with the last sentence of Remark \ref{rem.geodlum}, on the geometric irrelevance of additional parameter $(Z_0,T_0)$. Indeed, since by Equations $(10)$ we have ${\ds  \dot t_s \, =\, e^{-\sqrt{2}\,\omega\,x_s}\, b_s - a_s\,}$, we deduce from Corollary \ref{cor.gamma} that \  ${\ds X_s -\int_0^s (\dot z_u-\ell\,a_u) du\,}$ converges almost surely, and then, by (the proof of) Proposition \ref{pro.az''}, that so does also  
$$ X_s + \sigma \sqrt{1-\ell^2}\int_0^s \int_u^\ii e^{-\sigma^2 (v-u) - \sigma (w_{v} -w_u)} \, dw''_v\,du \, . $$ 
Now, as (setting $\,w^{(n)}_{\cdot} := w_{n+\cdot}-w_n\,$, for any fixed $n\in\N$) 
$$ \int_n^{n+1}\!\!\!\int_u^\ii\! e^{-\sigma^2 (v-u) - \sigma (w_{v} -w_u)}\,  dw''_v\,du\, =\, W''\bigg[\! \int_0^\ii \!\bigg[\! \int_0^{\min\{1,v\}}\! e^{\sigma^2 u + \sigma w^{(n)}_{u}} du \bigg]^2 e^{-2\sigma^2 v - 2 \sigma w^{(n)}_{v}} dv\bigg]  $$ 
has a constant law, we see that $(X_s)$ indeed cannot converge in probability. 
\if{ Note that however, by the same reason, 
$$\int_0^s \int_u^\ii e^{-\sigma^2 (v-u) - \sigma (w_{v} -w_u)} \, dw''_v\,du \,\equiv\, W''\bigg[\! \int_0^\ii \!\bigg[\! \int_0^{\min\{s,v\}}\! e^{\sigma^2 u + \sigma w_{u}} du \bigg]^2 e^{-2\sigma^2 v - 2 \sigma w_{v}} dv\bigg] $$ converges in law, proving that $(X_s)$ also does.  FAUX... }\fi 
}\erem 

\if{   Let us exhibit finally a fourth asymptotic random variable for the relativistic diffusion. 
FAUX... 
\bpro \label{pro.Z} \   The process \  ${\ds Z_s := \, z_s + \ell\, t_s - {2\,\ell\over \omega}\, {\rm Arctg}\left[ { 2-\rr\, e^{-\sqrt{2}\,\omega\,x_s} \over \sqrt{2(1-\ell^2)}\, + \omega\,\rr\, Y - \omega\,\rr\,y_s }\right]\,}$ \parn 
converges almost surely, as $\,s\to\ii$, toward some real random variable $\,Z\,$. 
\epro 
\ub{Proof} \quad  We have : 
$$ \dot Z_s = \dot z_s +\ell\, \dot t_s - 2\,\rr\,\ell\, { \sqrt{2}\,e^{-\sqrt{2}\,\omega\,x_s}\Big( \sqrt{2(1-\ell^2)}\, + \omega\,\rr\, Y - \omega\,\rr\,y_s \Big) \dot x_s +  \Big(2-\rr\, e^{-\sqrt{2}\,\omega\,x_s}\Big) \dot y_s \over  \Big(2-\rr\, e^{-\sqrt{2}\,\omega\,x_s}\Big)^2+ \Big( \sqrt{2(1-\ell^2)}\, + \omega\,\rr\, Y - \omega\,\rr\,y_s \Big)^2 } \, .  $$ 
Now, we know from Sections \ref{SubDiff2} and \ref{SubDiff4} that : 
$$ e^{-\sqrt{2}\,\omega\,x_s} = {\ts{2\over\rr}} - {\ts{\sqrt{2 (1-\ell^2)} \over \rr}}\,  \sin \gamma_s + o\Big(\! e^{-{2\over 3}\sigma^2 s}\!\Big) \; , \quad {\dot x_s/a_s} = \sqrt{1-\ell^2}\,\cos\gamma_s + o\Big(\! e^{-{2\over 3}\sigma^2 s}\!\Big)  , $$ 
from Equations $(10)$ that \quad ${\ds  \dot t_s \, =\, e^{-\sqrt{2}\,\omega\,x_s}\, b_s - a_s\, , } $ \quad $\dot y_s = e^{-\sqrt{2}\,\omega\,x_s} ( 2\,a_s- e^{-\sqrt{2}\,\omega\,x_s}b_s)$,  \parn 
and from Proposition \ref{pro.Y} that : \quad 
${\ds  y_s\, = \,Y_s - {\ts{\sqrt{2}\over \omega\, \rr}}\,{\dot x_s\over a_s} \, 
=\, Y_s - {\ts{\sqrt{2(1-\ell^2)}\over \omega\, \rr}}\, \cos\gamma_s + o\Big(\! e^{-{2\over 3}\sigma^2 s}\!\Big) . }$ \par \smallskip \noindent 
Hence we get : 
$$ {\dot Z_s\over a_s} = {\dot z_s\over a_s} +\ell\, e^{-\sqrt{2}\,\omega\,x_s} {b_s\over a_s} - \ell - 2\,\rr\,\ell\,e^{-\sqrt{2}\,\omega\,x_s}\, F_s\, , $$
with (using Remark \ref{rem.gamma}, and denoting by $\,o_s\,$ any exponentially fast vanishing remainder) : 
$$ F_s\, =\, { \sqrt{2}  \Big( {\ss\sqrt{2(1-\ell^2)}} + \omega \rr\, (Y - y_s) \Big) ({\ss\sqrt{1-\ell^2}}\,\cos\gamma_s + o_s) +  \Big(2-\rr\, e^{-\sqrt{2}\,\omega\,x_s}\Big) ( 2 - e^{-\sqrt{2}\,\omega\,x_s}{b_s\over a_s}) \over   \Big( \sqrt{2(1-\ell^2)}\,\sin\gamma_s + o_s \Big)^2+ \Big( \omega\,\rr \, (Y-Y_s) + \sqrt{2(1-\ell^2)}\,(1+\cos\gamma_s + o_s) \Big)^2  } $$ 
$$ =\, { 2\,(1-\ell^2)(1+\cos\gamma_s) \cos\gamma_s + {\ss\sqrt{2(1-\ell^2)}} \,\omega \rr\,(Y - Y_s) + o_s +  \Big({\ss\sqrt{2(1-\ell^2)}} \,\sin\gamma_s + o_s   \Big)^2 \over  4\,(1-\ell^2)(1+\cos\gamma_s) + 2\,{\ss\sqrt{2(1-\ell^2)}} \,\omega \rr\,(Y - Y_s) + \omega^2 \rr^2(Y - Y_s)^2 + o_s  } $$ 
$$ =\, { 2\,(1-\ell^2)(1+\cos\gamma_s) + {\ss\sqrt{2(1-\ell^2)}} \,\omega \rr\,(Y - Y_s) + o_s  \over  4\,(1-\ell^2)(1+\cos\gamma_s) + 2\,{\ss\sqrt{2(1-\ell^2)}} \,\omega \rr\,(Y - Y_s) + \omega^2 \rr^2(Y - Y_s)^2 + o_s  } $$ 
Now, by Proposition \ref{pro.Y} we have : \  $ (Y - Y_s) = o(a_s^{-4/9})$ 
$$ =\, { 2\, (1-\ell^2) (1+\cos\gamma_s)\, \cos\gamma_s+ \e_s +  \Big( 2- e^{\sqrt{2}\, \omega\,x_s} {b_s/a_s}\Big)\sqrt{2(1-\ell^2)}\, \sin\gamma_s \over  4\, (1-\ell^2) (1+\cos\gamma_s) + \e'_s   } $$ 
$$ =\,  { 2\, (1-\ell^2) (1+\cos\gamma_s)\, \cos\gamma_s +  2\, (1-\ell^2)\, \sin^2\gamma_s + \e_s'' \over  4\, (1-\ell^2) (1+\cos\gamma_s) + \e'_s   } \quad \hbox{ (by Remark \ref{rem.gamma})} $$ 
}\fi 
\medskip

   The theorem of the introduction (section \ref{sec.I}) is now established. Indeed, gathering successively Remark \ref{rem.lambda} and Proposition \ref{pro.irred}, Propositions \ref{pro.az} and \ref{pro.ell}, Corollary \ref{cor.x}, Proposition \ref{pro.Y} and Corollary \ref{cor.Y}, 
and Corollary \ref{cor.irred}, we get the following main result (for which $\,\sigma >0\,$ is necessary, due to the observation made after Definition \ref{def.convlum}, in Section \ref{sec.ngeod}).    

\bthe \label{the.compasym}  \   $(i)$ \  The relativistic diffusion is irreducible, on its phase space $\,\EE\,{\ss\sm}\,\EE_0\,$.  \par 
 $(ii)$ \  Almost surely, the relativistic diffusion path possesses a 3-dimensional asymptotic random variable $(\ell, \rr, Y)$, and converges to the light ray $\,B=(\ell,\rr, Y)\in \BB\,$, in the sense of Definition \ref{def.convlum}. Indeed, we have almost surely, as proper time $\,s\to\ii$ : 
$$ \dot z_s/a_s \lra \ell \,\in\, ]-1,0[\,\cup\,]0,1[\;;\quad b_s/a_s \lra  \rr\,\in \,]0,\ii[ \; ; \quad Y_s \lra  Y\,\in \R\, ;  $$ 
$$  \Big[ {\rr\over 2}\, e^{-\sqrt{2}\,\omega\,x_s} - 1\Big]^2 + \Big[{\omega\, \rr\over 2}\, (y_s - Y)\Big]^2 \lra\, \5(1-\ell^2)\, . $$ 

   $(iii)$ \  The asymptotic random variable $(\ell, \rr, Y)$ can be arbitrary close to any given $(\ell_0,\rr_0,y)\in\, ]-1,1[\times ]0,\ii[\times \R$, with positive probability. Hence,  the whole boundary (space of light rays) $\,\BB\,$ is the support of light rays the relativistic diffusion can converge to.   \ethe 

\brem \label{rem.dimB} \   {\rm  From the proofs of Lemma \ref{lem.ba}, Proposition \ref{pro.Y}, and Proposition \ref{pro.az}, we have the following representation of the asymptotic variable $\,B=(\ell,\rr,Y)$ : 
$$ \rr\, =\, {b_0\over a_0} + 2\sigma^2\! \int_0^\ii e^{\sqrt{2}\,\omega\,x_u}\,{du\over a_u^2} - \sigma^2\! \int_0^\ii {b_u\,du\over a_u^3} + \sigma\! \int_0^\ii\! a_u\1 \Big[ dM^b_u - {\ts {b_u\over a_u}}\, dM^a_u \Big] ;  $$ 
$$ Y\, =\, Y_0 - {\ts{2\sqrt{2}\,\sigma^2\over \omega}} \int_0^\ii  e^{2\sqrt{2}\,\omega\,x_u}\, {\dot x_u\over b_u^3}\, du + {\ts{\sqrt{2}\,\sigma\over \omega}} \int_0^\ii\! b_u\1 \Big[ dM^x_u - {\ts {\dot x_u\over b_u}}\, dM^b_u \Big] ;   $$ 
$$ \ell\, =\, {\dot z_0\over a_0} - \sigma^2 \int_{0}^\ii  {\dot z_u\over a_u^{3}}\, du + \sigma \, \int_0^\ii\! a_u\1 \Big[ dM^z_u - {\ts {\dot z_u\over a_u}}\, dM^a_u \Big] . $$ 
By Proposition \ref{lem.azn}, the law of the asymptotic variable $\,B\,$ has no atom, and by Theorem (\ref{the.compasym},$(iii)$), it is really three-dimensional. None of $\,\ell,\rr,Y\,$ is a function of the two others. 
}\erem 

   Recall that the random excitement of the relativistic diffusion is a standard three-dimensional Brownian motion. Therefore, Theorem \ref{the.compasym}, reinforced by Proposition \ref{pro.az''} and Remarks \ref{rem.ab'}, \ref{rem.lum}, \ref{rem.dimB}, incites to believe in the following.  \par \noindent 
\ub{\bf Conjecture} \  The tail $\sigma$-field and the invariant $\sigma$-field of the relativistic diffusion in G\"odel's universe are the $\sigma$-field generated by the asymptotic three-dimensional random variable $\, B=(\ell, \rr, Y)$ of Theorem \ref{the.compasym} (exhibited by Proposition \ref{pro.az}, Corollary \ref{cor.x}, and Proposition \ref{pro.Y}). 
\par 

\if{
\subsection{Scories} \label{sec.Scor} \indf 
   Fonctions testŽes : 
$$ {\dot x_s\over a_s} \; ; \quad  {\dot x_s\over a^2_s} \; ; \quad {\dot x_s^2\over a^2_s} \; ; \quad e^{-\sqrt{2}\,\omega\,x_s}\, {b_s\over a_s} \; ; \quad \Big[ e^{-\sqrt{2}\,\omega\,x_s}\, {b_s\over a_s} -2\Big]^2 \; ; \quad    e^{-\sqrt{2}\,\omega\,x_s}\, {b_s\, \dot x_s\over a_s^3}  \; ; \quad e^{-\sqrt{2}\,\omega\,x_s}\, {b_s\, \dot x_s\over a_s^2}  \; ; \quad $$ 
$$ e^{-\sqrt{2}\,\omega\,x_s}\, {\dot x_s\over a_s^n}  \; ; \quad \Big[ e^{-\sqrt{2}\,\omega\,x_s}\, {b_s\over a_s} -2\Big]\, {\dot x_s\over a_s^n}  \; ; \quad  {\dot x_s\over \sqrt{a_s}} \; ; \quad 
{\dot x_s^n\over a^p_s\,b_s^q} \; ; \quad  $$ 
ˆ tester ? 
$$ \Big[ e^{-\sqrt{2}\,\omega\,x_s}\, {b_s\over a_s} -2\Big]\, {\dot x_s^2\over a_s^n}  \; ; \quad   $$ 
Convergences obtenues : 
$$ \int_0^\ii\! \Big[ e^{-\sqrt{2}\,\omega\,x_u}\, {b_u\over a_u} -2\Big] e^{-\sqrt{2}\,\omega\,x_u}\,{b_u\over a_u}\, du\; ;  \;  \int_0^\ii\! \bigg[\Big[ e^{-\sqrt{2}\,\omega\,x_u}\, {\rr} -2\Big]^2+\ell^2-1-\e_u\bigg] e^{-\sqrt{2}\,\omega\,x_u}\,{b_u\over a_u}\, du\; .  $$ 
\blem \label{lem.b/a-rho} \  The integral \ ${\ds \int_0^\ii \Big[ {b_u\over a_u} -\rr\Big] e^{-\sqrt{2}\,\omega\,x_u}\,\dot x_u\,du\; }\,$ converges almost surely. 
\elem 
\ub{Proof} \quad  Integrating by parts, we have on one hand, for any $\,s\ge 0$ : 
$$ \sqrt{2}\,\omega \int_0^s \Big[ {b_u\over a_u} -\rr\Big] e^{-\sqrt{2}\,\omega\,x_u}\,\dot x_u\,du = 
\Big[ {b_0\over a_0} -\rr\Big] e^{-\sqrt{2}\,\omega\,x_0} - \Big[ {b_s\over a_s} -\rr\Big] e^{-\sqrt{2}\,\omega\,x_s} +  \int_0^s e^{-\sqrt{2}\,\omega\,x_u}\,d\Big[ {b_u\over a_u}\Big] \, , $$ 
with \ ${\ds \Big[ {b_s\over a_s} -\rr\Big] e^{-\sqrt{2}\,\omega\,x_s}}\;$ going to 0, \  and on the other hand 
$$ \int_0^s\! e^{-\sqrt{2}\,\omega\,x_u} d\Big[ {b_u\over a_u}\Big] = \sigma^2\!\int_0^s\! \Big[{2\over a_u^2} - {e^{-\sqrt{2}\,\omega\,x_u} b_u\over a_u^3}\Big] du  + \sigma\tilde W\bigg[\int_0^s\Big[ 4e^{-\sqrt{2}\,\omega\,x_u}\,{b_u\over a_u} -e^{-2\sqrt{2}\,\omega\,x_u} {b_u^2\over a_u^2} -2 \Big] {du \over a_u^2}\bigg] $$ 
converges almost surely, yielding the result. $\;\diamond$ 
\par \medskip  
}\fi 
\medskip

\brem \label{rem.conv} \   {\rm  The convergence of Theorem (\ref{the.compasym},$(ii)$), of the generic diffusion path $(\xi_s)$ to some light ray $\,\ol{\xi} = (\ol{t},\ol{x},\ol{y},\ol{z})$, in the sense of Definition \ref{def.convlum}, occurs in fact in some stronger sense. Indeed, by Remark \ref{rem.gamma}, Propositions \ref{pro.ab'} and \ref{pro.az}, and Lemma \ref{lem.a}, we have on one hand :
$$ e^{-\sqrt{2}\,\omega\,x_s} = {\ts{2\over \rr}} - {\ts{\sqrt{2(1-\ell^2)}\over \rr}} \,\sin\gamma_s + o(e^{-2\sigma^2s/3}) \;\hbox{ and } \quad y_s = Y - {\ts{\sqrt{2(1-\ell^2)}\over \omega\,\rr}} \,\cos\gamma_s + o(e^{-2\sigma^2s/3}) , $$ 
while by Remark \ref{rem.geodlum}, using the increasing diffeomorphism $\,\f=(\tau\mapsto\f_\tau)$, we have on the other hand : 
$$ e^{-\sqrt{2}\,\omega\,\ol{x}_{\f\1(\gamma_s/\omega)}} = {\ts{2\over \rr}} - {\ts{\sqrt{2(1-\ell^2)}\over \rr}} \,\sin\gamma_s \;\hbox{ and } \; \ol{y}_{\f\1(\gamma_s/\omega)} = Y - {\ts{\sqrt{2(1-\ell^2)}\over \omega\,\rr}} \,\cos\gamma_s \, . $$ 
Hence, we have in the $(x,y)$-plane a strong convergence, of the projection of the generic relativistic diffusion path to the projection of a lightlike geodesic, in the Skorohod topology : 
$$ \Big| x_s-\ol{x}_{\f\1(\gamma_s/\omega)}\Big| + \Big| y_s-\ol{y}_{\f\1(\gamma_s/\omega)}\Big|
= o(e^{-2\sigma^2s/3}) . $$ 
Otherwise, by Proposition \ref{pro.az''}, Corollaries \ref{cor.gamma} and \ref{cor.gamma'}, and Remark \ref{rem.geodlum}, we have : 
$$ z_s+\ell\, t_s = z_s+\ell t_0 + \ell \int_0^s\!(e^{-\sqrt{2}\,\omega\,x_u}b_u-a_u) du = {\ts{\ell\over \omega}}\, \gamma_s + \O(1) +  \int_0^s\!(\dot z_u-\ell\,a_u) du = {\ts{\ell\over \omega}}\, \gamma_s + o(|\gamma_s|^{1/3}) $$ 
$$ =\, \ol{z}_{\f\1(\gamma_s/\omega)}+\ell\,\ol{t}_{\f\1(\gamma_s/\omega)} + o(|\gamma_s|^{1/3}) =\, \Big(\ol{z}_{\f\1(\gamma_s/\omega)}+\ell\,\ol{t}_{\f\1(\gamma_s/\omega)}\Big)\Big[ 1 + o(e^{-\sigma^2 s/2})\Big] \lra \infty\, . $$ 
Hence, again in the Skorohod topology, and in the $(z,t)$-plane, the projection of the limiting light ray stands for an asymptotic direction for the projection of the generic relativistic diffusion path, but  there is no exactly asymptotic lightlike geodesic (only a parabolic branch occurs). 
}\erem 

\section{References} \label{Ref} 
\bigskip  

\vbox{ \noindent 
{\bf [B]} \ Bailleul I. \ \ {\it  Poisson boundary of a relativistic diffusion.}
\par \smallskip \hskip 24mm  Thesis of Orsay university, UniversitŽ Paris-Sud, France, 2006. } 
\medskip 

\vbox{ \noindent 
{\bf [D]} \ Dudley R.M. \ \ {\it Lorentz-invariant Markov processes in relativistic phase space.}
\par \smallskip \hskip 30mm  Arkiv f\"or Matematik 6, n$^o$ 14, 241-268, 1965. }
\medskip 

\if{ \vbox{ \noindent 
{\bf [D2]} \ Dudley R.M. \ \ {\it Recession of some relativistic Markov processes.}
\par \smallskip \hskip 33mm  Rocky Mountain J. Math. n$^o$ 4, 401-406, 1974. }
\medskip }\fi 

\vbox{  \noindent 
{\bf [F-LJ]} \ Franchi J. ,  Le Jan Y. \quad \hspace{-2mm}  {\it Relativistic Diffusions and Schwarzschild Geometry.} 
\par \hskip 54mm  ArXiv, \  http://arxiv.org/abs/math.PR/0410485, 2005.
\par \hskip 54mm  Comm. Pure Appl. Math., vol. LX, n$^o$ 2, 187-251, 2007. }
\medskip 

\vbox{  \noindent 
{\bf [G1]} \  G\"odel  K. \quad {\it An example of a new type of cosmological solution of Einstein's field  \par  \hskip 27mm equations of gravitation.} \quad  Rev. Mod. Phys. vol. 21, p. 447-450, 1949. } 
\medskip 

\vbox{  \noindent 
{\bf [G2]} \  G\"odel  K. \quad {\it Rotating universes in general relativity theory.} \par  \hskip 27mm   Proc.  Int. Congress Math., Cambridge, Mass., 1950, vol. 1,  p. 175-181. \par  \hskip 27mm
Amer. Math. Soc., Providence, R. I., 1952. } 
\medskip  

\vbox{ \noindent 
{\bf [H-E]} \ Hawking S.W. ,  Ellis G.F.R.\quad {\it The large-scale structure of space-time.}
\par \hskip 64mm Cambrige University Press, 1973. }
\medskip 

\vbox{  \noindent 
{\bf [I-W]} \ Ikeda N. ,  Watanabe S. \quad {\it Stochastic differential equations and diffusion
processes.}
\par \hskip 56mm North-Holland Kodansha, 1981. }
\medskip 

\vbox{ \noindent 
{\bf [K]} \   Kundt W. \quad {\it Tr\"agheitsbahnen in einem von G\"odel angegebenen kosmologischen Modell.} \par \hskip 26mm  Zeitschrift f\"ur Physik, vol. 145, p. 611-620, 1956.}
\medskip 

\vbox{ \noindent 
{\bf [M1]} \   Malament D. \  {\it Minimal acceleration requirement for time travel in G\"odel space-time.} \par \hskip 32mm  J. Math. Phys., vol. 26, n$^o\,$4, p. 774-777, 1985.}
\medskip 

\vbox{ \noindent 
{\bf [M2]} \   Malament D. \  {\it A note about closed timelike curves in G\"odel space-time.} \par \hskip 32mm  J. Math. Phys., vol. 28, n$^o\,$10, p. 2427-2430, 1987.}
\medskip 

\if{ 
\vbox{  \noindent 
{\bf [L-L]} \ Landau L. ,  Lifchitz E. \quad {\it Physique thŽorique, tome II : ThŽorie des champs.}
\par \hskip 54mm \'Editions MIR de Moscou, 1970. }
\medskip 

\vbox{ \noindent 
{\bf [M-T-W]} \ Misner C.W , Thorne K.S. , Wheeler J.A. \qquad    {\it Gravitation.} \par
\smallskip \hskip 15mm   W.H. Freeman and Company, New York, 1973. }
\medskip 

\vbox{  \noindent 
{\bf [P-R]} \ Penrose R. , Rindler W.  \quad {\it Spinors and space-time.}  \par 
\hskip 57mm  Cambridge University Press, 1986. }
\medskip 

\vbox{  \noindent 
{\bf [R]} \ Rindler W.  \quad {\it Relativity.}  \hskip 5mm  Oxford University Press, 2004. }
\medskip  

\vbox{  \noindent 
{\bf [W]} \ Wald R. M.  \quad {\it General relativity.}  \hskip 5mm  University of Chicago Press, 1984. }
\medskip  
 }\fi 
 
 \centerline{--------------------------------------------------------------------------------------------}
\medskip 
\noindent 
  \ub{A.M.S. Classification} : \ Primary 58J65, secondary 53C50, 60J65, 83C15, 83C10. \par 
\medskip \noindent 
  \ub{Key Words} : \  Relativistic diffusion, Brownian motion, Stochastic flow, General relativity, 
Lorentz manifold, G\"odel's universe, Light rays, Asymptotic behaviour. 

\medskip 
\centerline{--------------------------------------------------------------------------------------------}
\bigskip 

\vbox{  \noindent 
Jacques FRANCHI : \quad Universit\'e Louis Pasteur, I.R.M.A., 7 rue Ren\'e Descartes, \parn 
\hskip 43mm 67084 Strasbourg cedex. FRANCE. \quad \parn 
franchi@math.u-strasbg.fr 
}
\medskip 
\centerline{--------------------------------------------------------------------------------------------}

\end{document}